\documentclass[reqno]{siamart171218}
\usepackage{CJKutf8}
\usepackage[margin=0.86in]{geometry}

\usepackage{amssymb}
\usepackage{subfigure}
\usepackage{multirow}  
\usepackage{tabularx}

\usepackage{lipsum}
\usepackage{amsfonts}
\usepackage{graphicx}
\usepackage{epstopdf}
\usepackage{algorithmic}
\ifpdf
  \DeclareGraphicsExtensions{.eps,.pdf,.png,.jpg}
\else
  \DeclareGraphicsExtensions{.eps}
\fi

\newcommand{\mb}{\boldsymbol}
\newcommand{\xs}{\mathbb}

\newsiamremark{remark}{Remark}
\newsiamremark{hypothesis}{Hypothesis}
\crefname{hypothesis}{Hypothesis}{Hypotheses}
\newsiamthm{claim}{Claim}

\headers{AONN for parametric optimal control problems}{Pengfei Yin, Guangqiang xiao,  Kejun Tang and Chao Yang}

\title{AONN: An adjoint-oriented neural network method for all-at-once solutions of parametric optimal control problems\thanks{Submitted to the editors DATE.
}}

\author{Pengfei Yin\thanks{School of Mathematical Sciences, Peking University 
  (\email{pengfeiyin@pku.edu.cn}).}
\and Guangqiang Xiao\thanks{Academy for Advanced Interdisciplinary Studies, Peking University
  (\email{2001213226@stu.pku.edu.cn}).}
\and Kejun Tang\thanks{Peng Cheng Laboratory, Shenzhen, China
  (\email{tangkj@pcl.ac.cn}).}  
\and Chao Yang\thanks{School of Mathematical Sciences, Peking University}
  (\email{chao\_yang@pku.edu.cn}).}

\usepackage{amsopn}

\makeatletter
\newcommand*{\addFileDependency}[1]{
  \typeout{(#1)}
  \@addtofilelist{#1}
  \IfFileExists{#1}{}{\typeout{No file #1.}}
}
\makeatother

\newcommand*{\myexternaldocument}[1]{%
    \externaldocument{#1}%
    \addFileDependency{#1.tex}%
    \addFileDependency{#1.aux}%
}

\ifpdf
\hypersetup{
  pdftitle={
   AONN: An adjoint-oriented neural network method for all-at-once solutions of parametric optimal control problems},
  pdfauthor={Pengfei Yin, Guangqiang Xiao, Kejun Tang, Chao Yang}
}
\fi

\myexternaldocument{ex_supplement}

\usepackage{lineno}

\begin{document}

\maketitle

\begin{abstract}
Parametric optimal control problems governed by partial differential equations (PDEs) are widely found in scientific and engineering applications. Traditional grid-based numerical methods for such problems generally require repeated solutions of PDEs with different parameter settings, which is computationally prohibitive especially for problems with high-dimensional parameter spaces. 
Although recently proposed neural network methods make it possible to obtain the optimal solutions simultaneously for different parameters, challenges still remain when dealing with problems with complex constraints.
In this paper, we propose AONN, an adjoint-oriented neural network method, to overcome the limitations of existing approaches in solving parametric optimal control problems. In AONN, the neural networks are served as parametric surrogate models for the control, adjoint and state functions to get the optimal solutions all at once. 
In order to reduce the training difficulty and handle complex constraints, we introduce an iterative training framework inspired by the classical direct-adjoint looping (DAL) method so that penalty terms arising from the Karush–Kuhn–Tucker (KKT) system can be avoided.
Once the training is done, parameter-specific optimal solutions can be quickly computed through the forward propagation of the neural networks, which may be further used for analyzing the parametric properties of the optimal solutions. The validity and efficiency of AONN is demonstrated through a series of numerical experiments with problems involving various types of parameters.
\end{abstract}


\section{Introduction}
\label{sec:intro}
Optimal control modeling has been playing an important role in a wide range of applications, such as aeronautics \cite{xu2021machine}, mechanical engineering \cite{troltzsch2010optimal}, haemodynamics \cite{rozza2012reduction}, microelectronics \cite{pinnau2017quasi}, reservoir simulations \cite{yang2016nonlinear}, and environmental sciences \cite{strazzullo2018model}.
Particularly, to solve a PDE-constrained optimal control problem, one needs to find an optimal control function that can minimize a given cost functional for systems governed by partial differential equations (PDEs). 
Popular approaches for solving PDE-constrained optimal control problems include the direct-adjoint looping (DAL) method \cite{mitter1971optimal,jameson1988aerodynamic} that iteratively solves the adjoint systems, the Newton conjugate gradient method \cite{sternberg2010memory} that exploits the Hessian information, the semismooth Newton method \cite{yang2017nonlinearly, de2007optimal} that includes control and state constraints, and the alternating direction method of multipliers \cite{chen2019fe} designed for convex optimization.
In practice, the cost functionals and PDE systems often entail different configurations of physical or geometrical parameters, leading to parametric optimal control modeling. These parameters usually arise from certain desired profiles such as material properties, boundary conditions, control constraints, and computational domains \cite{karcher2016certified,karcher2018certified,strazzullo2018model, negri2015reduced,stadler2009elliptic,pinnau2017quasi,milani2008reduced}.

Most of the aforementioned methods cannot be directly applied to parametric optimal control problems.
The main reason is that, in addition to the already costly process of solving the PDEs involved in the optimal control modeling, the presence of parameters introduces extra prominent complexity, making the parametric optimal control problems much more challenging than the nonparametric ones \cite{hesthaven2016certified}. 
An efficient method for solving parametric optimal control problems is the reduced order model (ROM) \cite{quarteroni2015reduced, rozza2012reduction, strazzullo2018model, negri2013reduced}, which relies on surrogate models for the parametric model order reduction, and can provide both efficient and stable approximations if the solutions lie on a low-dimensional subspace \cite{hesthaven2016certified}.
However, because of the coupling of the spatial domain and the parametric domain,
the discretization in ROM still suffers from the curse of dimensionality, thus is unable to obtain all-at-once solutions to parametric optimal control problems \cite{karcher2016certified, rozza2012reduction, quarteroni2015reduced}.

Numerical methods based on deep learning have been receiving increasingly more attentions in solving PDEs \cite{raissi2018hidden,raissi2019physics, weinan2018deep,han2018solving, zhu2019physics,sheng2021pfnn,sheng2022pfnn}.
Recently, several successes have been made in solving PDE-constrained optimal control problems with deep-learning-based approaches.
For example, a physics-informed neural network (PINN) method is designed to solve optimal control problems by adding the cost functional to the standard PINN loss \cite{mowlavi2022optimal, lu2021physics}. 
Meanwhile, deep-learning-based surrogate models \cite{xu2021machine,lye2021iterative} and operator learning methods \cite{wang2021fast,hwang2022solving} are proposed to achieve fast inference for the optimal control solution without intensive computations. 
Although these methods are successful for solving optimal control problems, few of which can be directly applied in parametric optimal control modeling. In a recent work \cite{demo2021extended}, an extended PINN is proposed to augment neural network inputs with parameters, so that the Karush–Kuhn–Tucker (KKT) conditions and neural networks can be combined. 
In this way, the optimal solution with a continuous range of parameters could be obtained for parametric optimal control problems with simple constraints.
However, it is difficult for this method to generalize to solve more complex parametric optimal control problems, especially when the control function has additional inequality constraints \cite{ali2020reduced, bader2016certified}. 
In such scenarios, too many penalty terms have to be introduced into the loss function to fit the complex KKT system, which is very hard to optimize \cite{krishnapriyan2021characterizing}. 
A more detailed discussion of aforementioned methods can be found in Section~\ref{sec:analysis}.

To tackle the challenges in solving parametric optimal control problems and avoid the curse of dimensionality, we propose AONN, an adjoint-oriented neural network method that combines the advantages of both the classic DAL method and the deep learning technique.
In AONN, we construct three neural networks with augmented parameter inputs,  and integrate them into the framework of the DAL method to get an all-at-once approximation of the control function, the adjoint function, and the state function, respectively.
On the one hand, neural networks enable the classic DAL framework to solve parametric problems simultaneously with the aid of random sampling rather than the discretization of the coupled spatial domain and parametric domain. 
On the other hand, unlike the PINN-based penalty methods \cite{raissi2019physics,lu2021physics,mowlavi2022optimal,demo2021extended}, the introduction of DAL avoids directly solving the complex KKT system with various penalty terms.
Numerical results will show that, AONN can obtain high precision solutions to a series of  parametric optimal control problems.

The remainder of the paper is organized as follows. In Section~\ref{sec:problem}, the problem setting is introduced. After that, we will present the AONN framework in Section~\ref{sec:method}. Some further comparisons between AONN and several recently proposed methods are discussed in Section~\ref{sec:analysis}. Then, numerical results are presented in Section~\ref{sec:results} to demonstrate the efficiency of the proposed AONN method. The paper is concluded in Section~\ref{sec:conclusions}.
\section{Problem setup}
\label{sec:problem}
Let $\boldsymbol{\mu} \in \mathcal{P} \subset \xs{R}^D$ denote a vector that collects a finite number of parameters. Let $\Omega(\boldsymbol{\mu}) \subset \xs{R}^d$ be a spatial domain depending on $\boldsymbol{\mu}$, that is bounded, connected and with boundary $\partial \Omega(\boldsymbol{\mu})$, and $\mathbf{x} \in \Omega(\boldsymbol{\mu})$ denote a spatial variable.
Consider the following parametric optimal control problem
 \begin{linenomath*}\begin{equation}\label{OCPmu}
\mathrm{OCP}(\boldsymbol{\mu}):\quad
    \left\{\begin{aligned}
    &\min _{(y(\mathbf{x}, \boldsymbol{\mu}), u(\mathbf{x}, \boldsymbol{\mu})) \in Y \times U} J(y(\mathbf{x}, \boldsymbol{\mu}), u(\mathbf{x}, \boldsymbol{\mu}) ; \boldsymbol{\mu}),\\
    &\text { s.t. } \
    \mathbf{F}(y(\mathbf{x},  \boldsymbol{\mu}),u(\mathbf{x}, \boldsymbol{\mu});\boldsymbol{\mu}) = 0 \ \text{ in } \Omega(\boldsymbol{\mu}), 
    \text{ and } \ u(\mathbf{x}, \boldsymbol{\mu})\in U_{ad}(\boldsymbol{\mu}),
    \end{aligned}\right.
\end{equation}\end{linenomath*}
where $J: Y \times U \times \mathcal{P} \mapsto \mathbb{R}$ is a parameter-dependent objective functional, $Y$ and $U$ are two proper function spaces defined on $\Omega(\boldsymbol{\mu})$, with $y \in Y$ being the state function and $u \in U$ the control function, respectively. Both $y$ and $u$ are dependent on $\mathbf{x}$ and $\boldsymbol{\mu}$. To simplify the notation, we denote $y(\boldsymbol{\mu}) = y(\mathbf{x}, \boldsymbol{\mu})$ and $u(\boldsymbol{\mu}) = u(\mathbf{x}, \boldsymbol{\mu})$. 
In $\mathrm{OCP}(\boldsymbol{\mu})$ \eqref{OCPmu}, $\mathbf{F}$ represents the governing equation, such as, in our case, parameter-dependent PDEs, including the partial differential operator $\mathbf{F}_I$ and the boundary operator $\mathbf{F}_B$ (see Section~\ref{sec:results} for examples).
The admissible set $U_{ad}(\boldsymbol{\mu})$ is a parameter-dependent bounded closed convex subset of $U$, which provides an additional inequality constraint for $u$, e.g., the box constraint $U_{ad}(\boldsymbol{\mu})=\{u(\boldsymbol{\mu})\in U : u_a(\boldsymbol{\mu})\leq u(\boldsymbol{\mu}) \leq u_b(\boldsymbol{\mu})\}$.

Since the $\mathrm{OCP}(\boldsymbol{\mu})$ \eqref{OCPmu} is a constrained minimization problem, the necessary condition for the minimizer $(y^{*}(\boldsymbol{\mu}), u^*(\boldsymbol{\mu}))$ of \eqref{OCPmu} is the following KKT system \cite{troltzsch2010optimal,de2015numerical,hinze2008optimization}: 
 \begin{linenomath*}\begin{equation}\label{KKT}
    \left\{\begin{aligned}
     &J_y(y^{*}( \boldsymbol{\mu}),u^{*}( \boldsymbol{\mu}); \boldsymbol{\mu}) -  \mathbf{F}^{*}_{y}(y^{*}( \boldsymbol{\mu}),u^{*}(\boldsymbol{\mu});\boldsymbol{\mu})p^{*}(\boldsymbol{\mu}) =0, \\
    &\mathbf{F}(y^{*}( \boldsymbol{\mu}),u^{*}( \boldsymbol{\mu});\boldsymbol{\mu}) = 0, \\
     &(\mathrm{d}_u J(y^{*}(\boldsymbol{\mu}),u^{*}(\boldsymbol{\mu});\boldsymbol{\mu}),v(\boldsymbol{\mu})-u^{*}(\boldsymbol{\mu})) \geq 0, \ \forall v(\boldsymbol{\mu}) \in U_{ad}(\boldsymbol{\mu}),
    \end{aligned}\right.
\end{equation}\end{linenomath*}
where $p^{*}(\boldsymbol{\mu})$ is the adjoint function which is also known as the Lagrange multiplier, and $\mathbf{F}^{*}_{y}(y(\boldsymbol{\mu}),u(\boldsymbol{\mu});\boldsymbol{\mu})$ denotes the adjoint operator of $\mathbf{F}_{y}(y(\boldsymbol{\mu}),u(\boldsymbol{\mu});\boldsymbol{\mu})$. As $y(\boldsymbol{\mu})$ can always be uniquely determined by $u(\boldsymbol{\mu})$ through the state equation $\mathbf{F}$, the total derivative of $J$ with respect to $u$ in \eqref{KKT} can be formulated as 
 \begin{linenomath*}\begin{equation}\label{grad}
    \mathrm{d}_{u}J(y^{*}(\boldsymbol{\mu}),u^{*}(\boldsymbol{\mu});\boldsymbol{\mu}) = J_u(y^{*}(\boldsymbol{\mu}),u^{*}(\boldsymbol{\mu});\boldsymbol{\mu}) - \mathbf{F}_u^*(y^{*}(\boldsymbol{\mu}),u^{*}(\boldsymbol{\mu});\boldsymbol{\mu})p^{*}(\boldsymbol{\mu}).
\end{equation}\end{linenomath*}
The solution of OCP($\boldsymbol{\mu}$) satisfies the system \eqref{KKT}. So the key point is to solve this KKT system, based on which it is expected to find a minimizer for the OCP($\boldsymbol{\mu}$). In general, it is not a trivial task to solve \eqref{KKT} directly, and solving the parametric PDE involved in the KKT system poses additional computational challenges (e.g. the discretization of parametric spaces). In this work, we focus on the deep learning method to solve \eqref{KKT}. More specifically, we use three deep neural networks to approximate $y^{*}(\boldsymbol{\mu}), u^{*}(\boldsymbol{\mu})$ and $p^{*}(\boldsymbol{\mu})$ separately with an efficient training algorithm.
\section{Methodology}
\label{sec:method}
Let $\hat{y}\left(\mathbf{x}(\boldsymbol{\mu}); \boldsymbol{\theta}_{y}\right), \hat{u}\left(\mathbf{x}(\boldsymbol{\mu}); \boldsymbol{\theta}_{u}\right)$, and $\hat{p}\left(\mathbf{x}(\boldsymbol{\mu}); \boldsymbol{\theta}_{p}\right)$ be three independent deep neural networks parameterized with $\mb{\theta}_y, \mb{\theta}_u$ and $\mb{\theta}_p$ respectively. 
Here, $\mathbf{x}(\boldsymbol{\mu})$ is the augmented input of neural networks, which is given by
\begin{linenomath*}
\begin{equation*}
\mathbf{x}(\boldsymbol{\mu})=\left[\begin{array}{llllll}
x_{1}, & \ldots, & x_{d}, & \mu_{1}, & \ldots, & \mu_{D}
\end{array}\right].
\end{equation*}
\end{linenomath*}
We then use $\hat{y}\left(\mathbf{x}(\boldsymbol{\mu}); \boldsymbol{\theta}_{y}\right),  \hat{u}\left(\mathbf{x}(\boldsymbol{\mu}); \boldsymbol{\theta}_{u}\right)$, and $\hat{p}\left(\mathbf{x}(\boldsymbol{\mu}); \boldsymbol{\theta}_{p}\right)$ to approximate $y^{*}(\boldsymbol{\mu}), u^{*}(\boldsymbol{\mu})$ and $p^{*}(\boldsymbol{\mu})$ through minimizing three loss functions defined as
\begin{linenomath*}
\postdisplaypenalty=0
\begin{subequations}\label{eq_resloss}
\begin{align}
    \mathcal{L}_{s}(\boldsymbol{\theta}_y,\boldsymbol{\theta}_u) &= \left(\frac{1}{N}\sum_{i = 1}^{N} | r_{s}\left(\hat{y}(\mathbf{x}(\boldsymbol{\mu})_i;\mb{\theta}_y),\hat{u}(\mathbf{x}(\boldsymbol{\mu})_i;\mb{\theta}_u);\boldsymbol{\mu}_i\right) |^{2}\right)^{\frac{1}{2}},
    \label{lossstate}\\
    \mathcal{L}_{a}(\boldsymbol{\theta}_y,\boldsymbol{\theta}_u,\boldsymbol{\theta}_p) &= \left(\frac{1}{N}\sum_{i = 1}^{N} | r_{a}\left(\hat{y}(\mathbf{x}(\boldsymbol{\mu})_i;\mb{\theta}_y),\hat{u}(\mathbf{x}(\boldsymbol{\mu})_i;\mb{\theta}_u), \hat{p}(\mathbf{x}(\boldsymbol{\mu})_i;\mb{\theta}_p);\boldsymbol{\mu}_i\right) |^{2}\right)^{\frac{1}{2}},
    \label{lossadjoint}\\
    \mathcal{L}_{u}(\boldsymbol{\theta}_u, u_{\mathsf{step}}) &=  \left(\frac{1}{N}\sum_{i = 1}^{N} | \hat{u}(\mathbf{x}(\boldsymbol{\mu})_i;\mb{\theta}_u) - u_{\mathsf{step}}(\mathbf{x}(\boldsymbol{\mu})_i)|^{2}\right)^{\frac{1}{2}},
    \label{lossu}
\end{align}
\end{subequations}
\end{linenomath*}
where $\{\mathbf{x}(\boldsymbol{\mu})_i\}_{i=1}^N$ denote the collocation points.
The functionals $r_s$ and $r_a$ represent the residuals for the state equation and the adjoint equation induced by the KKT conditions, i.e.,
\begin{linenomath*}
\postdisplaypenalty=0
\begin{subequations}\label{res}
\begin{align}
r_{s}(y(\boldsymbol{\mu}),u(\boldsymbol{\mu});\boldsymbol{\mu}) &\triangleq \mathbf{F}(y(\boldsymbol{\mu}),u(\boldsymbol{\mu});\boldsymbol{\mu}), \label{res_s}\\
r_{a}(y(\boldsymbol{\mu}),u(\boldsymbol{\mu}),p(\boldsymbol{\mu});\boldsymbol{\mu}) &\triangleq J_y(y(\boldsymbol{\mu}),u(\boldsymbol{\mu});\boldsymbol{\mu}) - \mathbf{F}^{*}_{y}(y(\boldsymbol{\mu}),u(\boldsymbol{\mu});\boldsymbol{\mu})p(\boldsymbol{\mu}), \label{res_a}
\end{align}
\end{subequations}
\end{linenomath*}
and $u_{\mathsf{step}}(\mathbf{x}(\boldsymbol{\mu}))$ is an intermediate variable during the update procedure of the control function for the third variational inequality in the KKT conditions \eqref{KKT}, which will be discussed in Section~\ref{sec_proj}. These three loss functions try to fit the KKT conditions by adjusting the parameters of the three neural networks $\hat{y}\left(\mathbf{x}(\boldsymbol{\mu}); \boldsymbol{\theta}_{y}\right),  \hat{u}\left(\mathbf{x}(\boldsymbol{\mu}); \boldsymbol{\theta}_{u}\right)$, and $\hat{p}\left(\mathbf{x}(\boldsymbol{\mu}); \boldsymbol{\theta}_{p}\right)$, and the training procedure is performed in a sequential way. The derivatives involved in the loss functions can be computed efficiently by automatic differentiation in deep learning libraries such as TensorFlow \cite{abadi2016tensorflow} or PyTorch \cite{paszke2017automatic}.

\subsection{Deep learning for parametric PDEs}
The efficient solution of parametric PDEs is crucial for parametric optimal control modeling because of extra parameters involved in the physical system \eqref{OCPmu}.
To deal with parametric PDEs, we augment the input space of the neural networks by taking the parameter $\boldsymbol{\mu}$ as additional inputs, along with the coordinates $\mathbf{x}$ to handle the parameter-dependent PDEs. In addition, the penalty-free techniques \cite{lagaris1998artificial,sheng2021pfnn} are employed to enforce boundary conditions in solving parametric PDEs.
Next, we illustrate how to apply penalty-free deep neural networks to solve the parametric state equation \eqref{res_s}, which can be directly generalized to the solution of the adjoint equation \eqref{res_a}.

The key point of the penalty-free method is to introduce two neural networks to approximate the solution, of which one neural network $\hat{y}_{B}$ is used to approximate the essential boundary conditions and the other $\hat{y}_{I}$ deals with the rest part of the computational domain. 
In this way, the training difficulties from the boundary conditions are eliminated, which improves the accuracy and robustness for complex geometries.
For problems with simple geometries, we can also construct an analytical expression for $\hat{y}_{B}$ to further reduce the training cost (see Section~\ref{sec:results} for examples).
The approximate solution of the state equation is constructed by
 \begin{linenomath*}\begin{equation}\label{eq-pfnn}
\hat{y}\left(\mathbf{x}(\boldsymbol{\mu}) ; \boldsymbol{\theta}_{y}\right) = \hat{y}_{B}(\mathbf{x}(\boldsymbol{\mu});\mb{\theta}_{y_B})+\ell(\mathbf{x}(\boldsymbol{\mu})) \hat{y}_{I}\left(\mathbf{x}(\boldsymbol{\mu}) ; \boldsymbol{\theta}_{y_I}\right),
\end{equation}\end{linenomath*}
where $\mb{\theta}_y = \{\mb{\theta}_{y_B}, \mb{\theta}_{y_I}\}$ collects all parameters of two sub-neural networks $\hat{y}_{B}$ and $\hat{y}_{I}$, and $\ell$ is a length factor function that builds the connection between $\hat{y}_{B}$ and $\hat{y}_{I}$, satisfying the following two conditions:
\begin{linenomath*}
\begin{equation*}
    \left\{\begin{array}{l}
    \ell(\mathbf{x}(\boldsymbol{\mu})) > 0, \quad \mbox{in}\;\Omega(\boldsymbol{\mu}),\\
    \ell(\mathbf{x}(\boldsymbol{\mu})) = 0 ,\quad \mbox{on} \;\partial \Omega(\boldsymbol{\mu}).
    \end{array}\right.\\
\end{equation*}
\end{linenomath*}
The details of constructing the length factor function $\ell$ can be found in ref.\cite{sheng2021pfnn}. With these settings, training $\hat{y}_{B}$ and $\hat{y}_{I}$ can be performed separately, i.e., one can first train $\hat{y}_{B}$, and then fix $\hat{y}_{B}$ to train $\hat{y}_{I}$.
For a fixed $u(\boldsymbol{\mu})$, we have 
\begin{linenomath*}
\begin{equation*}
\begin{aligned}
\mathbf{F}(\hat{y}\left(\mathbf{x}(\boldsymbol{\mu}) ; \boldsymbol{\theta}_{y}\right), u(\boldsymbol{\mu});\boldsymbol{\mu}) =
\left[\begin{array}{l}
\mathbf{F}_I(\hat{y}\left(\mathbf{x}(\boldsymbol{\mu}) ; \boldsymbol{\theta}_{y}\right),u(\boldsymbol{\mu});\boldsymbol{\mu}) \\
\mathbf{F}_B(\hat{y}_{B}(\mathbf{x}(\boldsymbol{\mu});\mb{\theta}_{y_B}),u(\boldsymbol{\mu});\boldsymbol{\mu})
\end{array}\right],
\end{aligned}
\end{equation*}
\end{linenomath*}
and the residual of the state equation can be rewritten as
\begin{linenomath*}
\begin{equation*}
r_{s}(\hat{y}\left(\mathbf{x}(\boldsymbol{\mu}) ; \boldsymbol{\theta}_{y}\right), u(\boldsymbol{\mu});\boldsymbol{\mu}) =  
\left[\begin{array}{l}
r_{s_I}(\hat{y}\left(\mathbf{x}(\boldsymbol{\mu}) ; \boldsymbol{\theta}_{y}\right),u(\boldsymbol{\mu});\boldsymbol{\mu}) \\
r_{s_B}(\hat{y}_{B}(\mathbf{x}(\boldsymbol{\mu});\mb{\theta}_{y_B}),u(\boldsymbol{\mu});\boldsymbol{\mu})
\end{array}\right].
\end{equation*}
\end{linenomath*}
We then sample a set $\{\mathbf{x}(\boldsymbol{\mu})_i\}_{i=1}^N$ of collocation points to optimize $\mb{\theta}_y$ through minimizing the state loss function \eqref{lossstate} if $\mb{\theta}_u$ is fixed.

For parametric problems, we take parameters as the additional inputs of neural networks. This approach is used to solve parametric forward problems \cite{khodayi2020varnet} and control problems \cite{sun2020surrogate}. A typical way for sampling training points is to separately sample data in $\Omega$ and $\mathcal{P}$ to get $\{\mathbf{x}_i\}$ and $\{\boldsymbol{\mu}_j\}$, and then compose  product data $\{(\mathbf{x}_i, \boldsymbol{\mu}_j)\}$ in $\Omega\times\mathcal{P}$. Rather than taken from each slice of $\Omega(\boldsymbol{\mu})$ for a fixed $\boldsymbol{\mu}$, collocation points are sampled in space $\Omega_{\mathcal{P}}$ in this work, where
\begin{linenomath*}
\begin{equation*}
    \Omega_{\mathcal{P}} = \{\mathbf{x}(\boldsymbol{\mu}):\mathbf{x}\in \Omega(\boldsymbol{\mu})\}
\end{equation*}
\end{linenomath*}
represents the joint spatio-parametric domain. The reason is that, for parametric geometry problems, as the spatial domain $\Omega(\boldsymbol{\mu})$ is parameter-dependent, the sampling space cannot be expressed as the Cartesian product of $\Omega$ and $\mathcal{P}$.

\subsection{Projection gradient descent} \label{sec_proj}
Due to the additional inequality constraints $ u \in U_{ad}(\boldsymbol{\mu})$ for the control function, the zero gradient condition $\mathrm{d}_{u}J(y^{*}(\boldsymbol{\mu}),u^{*}(\boldsymbol{\mu});\boldsymbol{\mu})=0$ cannot be directly applied to the optimal solution to get the update scheme for $u$. One way to resolve this issue is to introduce additional Lagrange multipliers with some slack variables to handle the inequality constraints. However, this will bring additional penalty terms that could affect the procedure of optimization \cite{chen2020multi}. Furthermore, the inequality constraints also lead to the non-smoothness of the control function, making it more difficult to capture the singularity by penalty methods \cite{lu2021physics, haghighat2021physics}. To avoid these issues, we here use a simple iterative method to handle the variational inequality without utilizing a Lagrange multiplier, where a projection gradient descent method is employed, based on which we can obtain the update scheme for $u$.  The projection operator onto the admissible set $U_{ad}(\boldsymbol{\mu})$ is defined as:
\begin{linenomath*}
\begin{equation*}
    \mathbf{P}_{U_{ad}(\boldsymbol{\mu})}(u(\boldsymbol{\mu})) = \arg\min_{v(\boldsymbol{\mu})\in U_{ad}(\boldsymbol{\mu})}\|u(\boldsymbol{\mu})-v(\boldsymbol{\mu})\|_{2},
\end{equation*}
\end{linenomath*}
which performs the projection of $u(\boldsymbol{\mu})$ onto the convex set $U_{ad}(\boldsymbol{\mu})$. In practice, the above projection is implemented in a finite dimensional vector space, i.e., $u(\boldsymbol{\mu})$ is discretized on a set of collocation points (e.g. grids on the domain $\Omega$). So it is straightforward to build this projection since $U_{ad}(\boldsymbol{\mu})$ is a convex set.
For example, if 
 \begin{linenomath*}\begin{equation}\label{box-constraint}
    U_{ad}(\boldsymbol{\mu}) = \{u\in U:u_a(\mathbf{x}(\boldsymbol{\mu}))\leq u(\mathbf{x}(\boldsymbol{\mu})) \leq u_b(\mathbf{x}(\boldsymbol{\mu})),\forall \mathbf{x} \in \Omega(\boldsymbol{\mu}) \}
\end{equation}\end{linenomath*}
provides a box constraint for $u$, where $u_a$ and $u_b$ are the lower bound function and the upper bound function respectively, both of which are dependent on $\boldsymbol{\mu}$, and $[u_1, \ldots, u_N]^{\mathsf{T}}$ represents the control function values at $N$ collocation points $\{\mathbf{x}(\boldsymbol{\mu})_{i}\}_{i=1}^{N}$ in $\Omega_{\mathcal{P}}$, then we can construct the projection $\mathbf{P}_{U_{ad}(\boldsymbol{\mu})}$ in an entry-wise way \cite{troltzsch2010optimal,hinze2008optimization}:
 \begin{linenomath*}\begin{equation}\label{eqn_proj}
		\mathbf{P}_{U_{ad}(\boldsymbol{\mu})}(u_i) = \begin{cases} u_a(\mathbf{x}(\boldsymbol{\mu})_{i}), \quad  \mathrm{if} \  u_i < u_a(\mathbf{x}(\boldsymbol{\mu})_{i}) ,\\
			u_i,  \quad \mathrm{if}  \  u_b(\mathbf{x}(\boldsymbol{\mu})_{i}) \geq u_i \geq u_a(\mathbf{x}(\boldsymbol{\mu})_{i}) , \\
			u_b(\mathbf{x}(\boldsymbol{\mu})_{i}), \quad  \mathrm{if} \ u_i > u_b(\mathbf{x}(\boldsymbol{\mu})_{i}).
		\end{cases} i = 1, \ldots, N.
\end{equation}\end{linenomath*}
The projection gradient step can be carried out according to the above formula, so as to obtain the update of the control function denoted by $u_{\mathsf{step}}$, which is
 \begin{linenomath*}\begin{equation} \label{eq_uproj}
 u_{\mathsf{step}}(\boldsymbol{\mu}) = \mathbf{P}_{U_{ad}(\boldsymbol{\mu})}\left(u(\boldsymbol{\mu}) - c \mathrm{d}_{u}J(y(\boldsymbol{\mu}),u(\boldsymbol{\mu});\boldsymbol{\mu})\right),
\end{equation}\end{linenomath*}
and the loss for updating the control function is naturally defined as in \eqref{lossu}, making an approximation of $\hat{u}$ obtained through minimizing \eqref{lossu}.

The optimal control function $u^*(\boldsymbol{\mu})$ satisfies the following variational property:
\begin{linenomath*}
\begin{equation*}
u^*(\boldsymbol{\mu}) - \mathbf{P}_{U_{ad}(\boldsymbol{\mu})}\left(u^*(\boldsymbol{\mu})-c \mathrm{d}_{u}J(y^*(\boldsymbol{\mu}),u^*(\boldsymbol{\mu});\boldsymbol{\mu})\right) = 0, \quad \forall c\geq 0.
\end{equation*}
\end{linenomath*}
Here $\mathrm{d}_{u}J$ is associated with the adjoint function $p^{*}(\boldsymbol{\mu})$ from total derivative expression \eqref{grad}, and thus we define the residual for the control function
 \begin{linenomath*}\begin{equation}\label{eq_rescontrol}
r_{v}(y(\boldsymbol{\mu}),u(\boldsymbol{\mu}),p(\boldsymbol{\mu}),c;\boldsymbol{\mu})\triangleq u(\boldsymbol{\mu}) - \mathbf{P}_{U_{ad}(\boldsymbol{\mu})}\left(u(\boldsymbol{\mu})-c \mathrm{d}_{u}J(y(\boldsymbol{\mu}),u(\boldsymbol{\mu});\boldsymbol{\mu})\right),
\end{equation}\end{linenomath*}
and its corresponding variational loss is defined as
 \begin{linenomath*}\begin{equation}\label{lossvar}
\mathcal{L}_{v}(\boldsymbol{\theta}_y,\boldsymbol{\theta}_u,\boldsymbol{\theta}_p,c) = \left(\frac{1}{N}\sum_{i = 1}^{N} | r_{v}\left(\hat{y}(\mathbf{x}(\boldsymbol{\mu})_i;\mb{\theta}_y),\hat{u}(\mathbf{x}(\boldsymbol{\mu})_i;\mb{\theta}_u), \hat{p}(\mathbf{x}(\boldsymbol{\mu})_i;\mb{\theta}_p),c;\boldsymbol{\mu}_i\right) |^{2}\right)^{\frac{1}{2}}.
\end{equation}\end{linenomath*}

The first two losses in \eqref{eq_resloss} together with \eqref{lossvar} reflect how well $y(\boldsymbol{\mu}), u(\boldsymbol{\mu})$ and $p(\boldsymbol{\mu})$ approximate the optimal solution governed by the KKT system \eqref{KKT}.
Note that $r_{v}$ and $\mathcal{L}_{v}$ are dependent on the constant $c$, which is actually the step size for gradient descent. For verification, the variational loss is constructed to verify the convergence of algorithm, and $c$ is often chosen as the last step size.

\subsection{AONN algorithm}
\begin{figure}[htbp]
\label{AONN-flow}
\centering
\includegraphics[height=7.5cm]{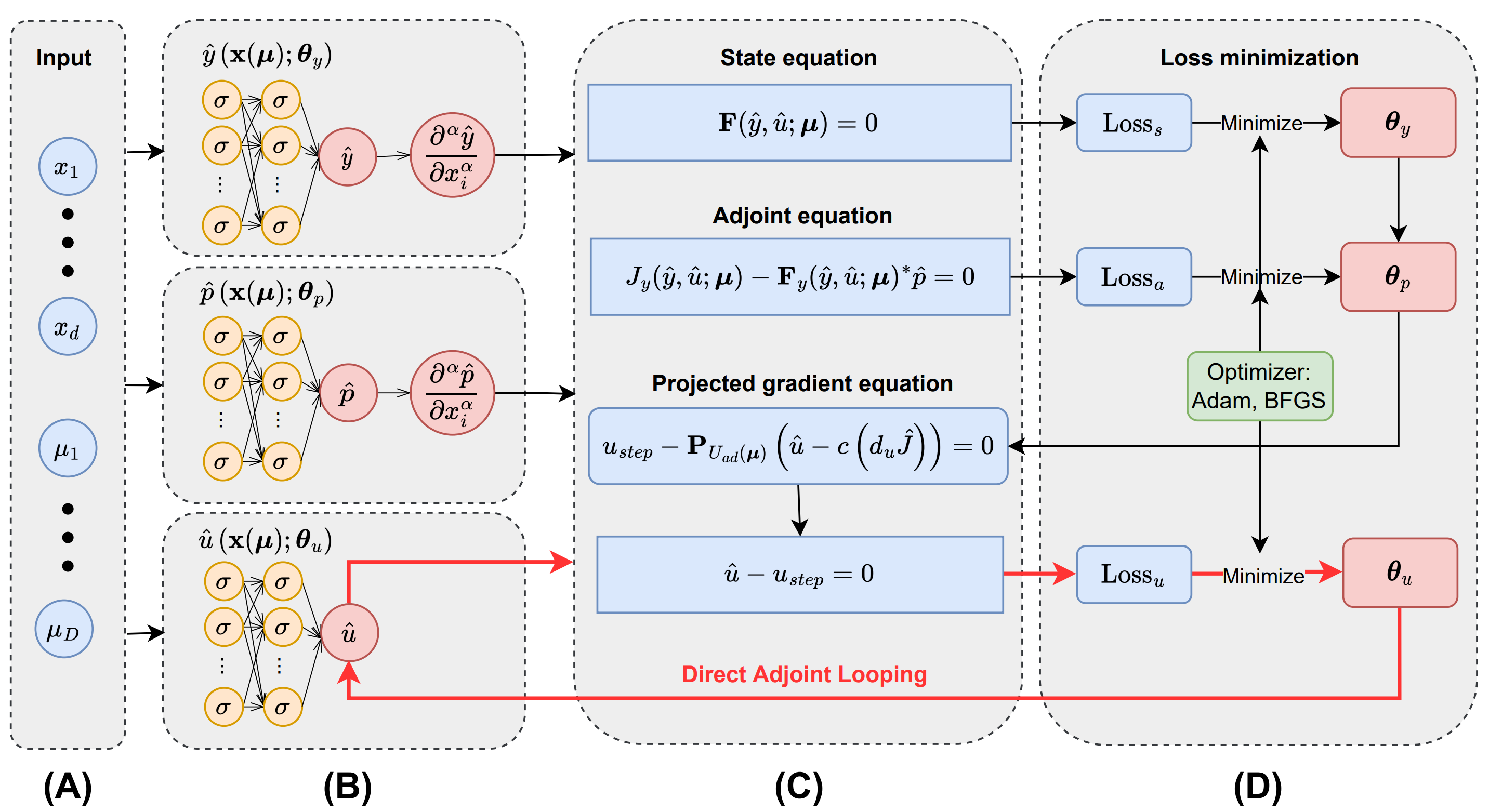}
\caption{$\text{\rm{The schematic of AONN for solving the parametric optimal control problems.}}$ $\rm{(A)}$ Spatial coordinates and parameters form the input of neural networks. $\rm{(B)}$ AONN consists of three separate neural networks $\hat{y}, \hat{p}, \hat{u}$ and return the approximation of state, adjoint and control respectively. $\rm{(C)}$ The state equation, the adjoint eqaution and the projected gradient equation are derived to formulate the corresponding loss functions. $\rm{(D)}$ The gradients in the state PDE and the adjoint PDE are computed via automatic differentiation \cite{paszke2017automatic}. $\hat{y},\hat{p},\hat{u}$ are then trained sequentially via the Adam \cite{kingma2014adam} or the BFGS optimizer.}
\end{figure}

Now putting all together, we are ready to present our algorithm. Our goal is to efficiently approximate the minimizer of \eqref{OCPmu} via adjoint-oriented neural networks (AONN). The overall training procedure of AONN consists of three steps: training $\hat{y}\left(\mathbf{x}(\boldsymbol{\mu}) ; \boldsymbol{\theta}_{y}\right)$, updating $\hat{p}\left(\mathbf{x}(\boldsymbol{\mu}) ; \boldsymbol{\theta}_{p}\right)$  and refining $\hat{u}\left(\mathbf{x}(\boldsymbol{\mu}) ; \boldsymbol{\theta}_{u}\right)$. 
The schematic of AONN for solving the parametric optimal control problems is shown in Figure~\ref{AONN-flow}.
The three neural networks $(\hat{y}\left(\mathbf{x}(\boldsymbol{\mu}) ; \boldsymbol{\theta}_{y}\right), \hat{p}\left(\mathbf{x}(\boldsymbol{\mu}) ; \boldsymbol{\theta}_{p}\right), \hat{u}\left(\mathbf{x}(\boldsymbol{\mu}) ; \boldsymbol{\theta}_{u}\right))$ with augmented parametric input, as illustrated in panels A and B, are optimized to iteratively minimizing the objective functional with respect to the corresponding variables once at a time. More specifically, according to the loss functions derived by the three equations shown in panel C, the training procedure is performed as in panel D. 

Starting with three initial neural networks $\hat{y}\left(\mathbf{x}(\boldsymbol{\mu}) ; \boldsymbol{\theta}^0_{y}\right), \hat{p}\left(\mathbf{x}(\boldsymbol{\mu}) ; \boldsymbol{\theta}^0_{p}\right)$, and $\hat{u}\left(\mathbf{x}(\boldsymbol{\mu}) ; \boldsymbol{\theta}^0_{u}\right)$, we train and obtain the state function $\hat{y}\left(\mathbf{x}(\boldsymbol{\mu}) ; \boldsymbol{\theta}^1_{y}\right)$ through minimizing $\mathcal{L}_{s}\left(\boldsymbol{\theta}_{y}, \boldsymbol{\theta}_{u}^0\right)$ (see \eqref{lossstate}), which is equivalent to solving the parameter-dependent state equation. With $\hat{y}\left(\mathbf{x}(\boldsymbol{\mu}) ; \boldsymbol{\theta}^1_{y}\right)$, we minimize the loss $\mathcal{L}_{a}\left(\boldsymbol{\theta}_{y}^{1}, \boldsymbol{\theta}_{u}^{0}, \boldsymbol{\theta}_{p}\right)$ (see \eqref{lossadjoint}) for the adjoint equation to get $\hat{p}\left(\mathbf{x} ; \boldsymbol{\theta}^1_{p}\right)$, corresponding to solving the parameter-dependent adjoint equation. To update the control function $\hat{u}\left(\mathbf{x}(\boldsymbol{\mu}) ; \boldsymbol{\theta}_{u}\right)$, $u^0_{\mathsf{step}}(\mathbf{x}(\boldsymbol{\mu}))$ is computed first by gradient descent followed by a projection step (see \eqref{eq_uproj}), and then $\hat{u}\left(\mathbf{x}(\boldsymbol{\mu}) ; \boldsymbol{\theta}^1_{u}\right)$ is obtained by minimizing $\mathcal{L}_{u}(\boldsymbol{\theta}_u, u^0_{\mathsf{step}}(\mathbf{x}(\boldsymbol{\mu})))$ (see \eqref{lossu}). Then another iteration starts using $\mb{\theta}_y^1, \mb{\theta}_p^1, \mb{\theta}_u^1$ as the initial parameters.  
In general, the iterative scheme is specified as follows:
\begin{linenomath*}
\begin{equation*}
\begin{aligned}
  \text{training } \hat{y}: &&\boldsymbol{\theta}_{y}^{k} &= \arg \min_{\boldsymbol{\theta}_{y}} \mathcal{L}_{s}\left(\boldsymbol{\theta}_{y}, \boldsymbol{\theta}_{u}^{k-1}\right),\\
  \text{updating } \hat{p}: &&\boldsymbol{\theta}_{p}^{k} &= \arg \min_{\boldsymbol{\theta}_{p}} \mathcal{L}_{a}\left(\boldsymbol{\theta}_{y}^{k}, \boldsymbol{\theta}_{u}^{k-1}, \boldsymbol{\theta}_{p}\right),\\
  \text{refining } \hat{u}: &&\boldsymbol{\theta}_{u}^{k} &= \arg \min_{\boldsymbol{\theta}_{u}} \mathcal{L}_{u}\left(\boldsymbol{\theta}_u, u^{k-1}_{\mathsf{step}}\right),
\end{aligned}
\end{equation*}
\end{linenomath*}
where 
 \begin{linenomath*}\begin{equation}
u^{k-1}_{\mathsf{step}}(\mathbf{x}(\boldsymbol{\mu})) = 
\mathbf{P}_{U_{ad}(\boldsymbol{\mu})}\left(\hat{u}(\mathbf{x}(\boldsymbol{\mu});\boldsymbol{\theta}_{u}^{k-1})-c^k
	    \mathrm{d}_{u}J(\hat{y}(\mathbf{x}(\boldsymbol{\mu});\boldsymbol{\theta}_{y}^{k}),\hat{u}(\mathbf{x}(\boldsymbol{\mu});\boldsymbol{\theta}_{u}^{k-1});\boldsymbol{\mu})\right),
\label{ustep1}
\end{equation}\end{linenomath*}
and
 \begin{linenomath*}\begin{equation}\label{ustep2}
\begin{aligned}
\mathrm{d}_{u}J\left(\hat{y}(\mathbf{x}(\boldsymbol{\mu});\boldsymbol{\theta}_{y}^{k}),\hat{u}(\mathbf{x}(\boldsymbol{\mu});\boldsymbol{\theta}_{u}^{k-1});\boldsymbol{\mu}\right) &= J_u\left(\hat{y}(\mathbf{x}(\boldsymbol{\mu});\boldsymbol{\theta}_{y}^{k}),\hat{u}(\mathbf{x}(\boldsymbol{\mu});\boldsymbol{\theta}_{u}^{k-1});\boldsymbol{\mu}\right) \\
&- \mathbf{F}_u^*\left(\hat{y}(\mathbf{x}(\boldsymbol{\mu});\boldsymbol{\theta}_{y}^{k}),\hat{u}(\mathbf{x}(\boldsymbol{\mu});\boldsymbol{\theta}_{u}^{k-1});\boldsymbol{\mu}\right)\hat{p}(\mathbf{x}(\boldsymbol{\mu});\boldsymbol{\theta}_{p}^{k}).
\end{aligned}
\end{equation}\end{linenomath*}

The iteration of AONN revolves around the refinement of $\hat{u}$ with the aid of $\hat{y}$ and $\hat{p}$, forming the direct-adjoint looping (DAL) ,which is indicated by red lines in Figure~\ref{AONN-flow}.
This procedure shares similarities to the classical DAL framework, but there is a crucial difference between AONN and DAL. That is, a reliable solution of $\mathrm{OCP}(\boldsymbol{\mu})$ for any parameter can be efficiently computed from the trained neural networks in our AONN framework, while DAL cannot achieve this. More details can be found in the discussions of Section~\ref{sec:analysis}.

The training process is summerized in Algorithm~\ref{alg_AONN}, where the loss function $\mathcal{L}_{v}(\boldsymbol{\theta}_y,\boldsymbol{\theta}_u,\boldsymbol{\theta}_p,c)$ (see \eqref{lossvar}) is used for the verification. In our practical implementation, we employ the step size decay technique with a decay factor $\gamma$ for robustness. The AONN method can be regarded as an inexact DAL to some extent since the state equation and the adjoint equation are not accurately solved but approximated with neural networks at each iteration. So the number of epochs is increased by $n_{\mathrm{aug}}$ compared with the previous step (on line 9 of Algorithm~\ref{alg_AONN}) to ensure the accuracy and convergence. It is worth noting that the training of  network $\hat{u}\left(\mathbf{x}(\boldsymbol{\mu}) ; \boldsymbol{\theta}^*_{u}\right)$ can be put after the while loop, if the collocation points are always fixed, since training the state function only uses the value of $u$ at the collocation points (the calculation of line $5$ in Algorithm~\ref{alg_AONN}).

\begin{algorithm}[H]
\caption{AONN for $\mathrm{OCP(\boldsymbol{\mu})}$}
\label{alg_AONN}
\begin{algorithmic}[1]
	\REQUIRE Initial $\boldsymbol{\theta}_y^0, \boldsymbol{\theta}_u^0, \boldsymbol{\theta}_p^0$,
	 collocation points$\{\mathbf{x}(\boldsymbol{\mu})_i\}_{i=1}^N$,
	 decay factor $\gamma \in (0,1]$, initial step size $c^0$, initial number of epochs $n^0$, positive integer $n_{\mathrm{aug}}$ and total iterations $N_{\mathrm{iter}}$.
	\STATE $k \longleftarrow 1$
	\WHILE{$k<N_{\mathrm{iter}}+1$}
	    \STATE $\boldsymbol{\theta}_{y}^{k} \longleftarrow \arg \min_{\boldsymbol{\theta}_{y}} \mathcal{L}_{s}\left(\boldsymbol{\theta}_{y}, \boldsymbol{\theta}_{u}^{k-1}\right)$: Train network $\hat{y}\left(\mathbf{x}(\boldsymbol{\mu}) ; \boldsymbol{\theta}_{y}\right)$ with initialization $\boldsymbol{\theta}_{y}^{k-1}$ for $n^k$ epochs. 
	    
	    \STATE $\boldsymbol{\theta}_{p}^{k} \longleftarrow \arg \min _{\boldsymbol{\theta}_{p}} \mathcal{L}_{a}\left(\boldsymbol{\theta}_{y}^{k}, \boldsymbol{\theta}_{u}^{k-1}, \boldsymbol{\theta}_{p}\right)$: Train network $\hat{p}\left(\mathbf{x}(\boldsymbol{\mu}) ; \boldsymbol{\theta}_{p}\right)$ with initialization $\boldsymbol{\theta}_{p}^{k-1}$ for $n^k$ epochs. 
	    \STATE Compute $u_{\mathsf{step}}^{k-1}(\mathbf{x}(\boldsymbol{\mu}))$ by \eqref{ustep1} and \eqref{ustep2}.
        \STATE $\boldsymbol{\theta}_{u}^{k} \longleftarrow \arg \min _{\boldsymbol{\theta}_{u}} \mathcal{L}_{u}(\boldsymbol{\theta}_u, u^{k-1}_{\mathsf{step}})$: Train  network $\hat{u}\left(\mathbf{x}(\boldsymbol{\mu}) ; \boldsymbol{\theta}_{u}\right)$ with initialization $\boldsymbol{\theta}_{u}^{k-1}$ for $n^k$ epochs. 
		\STATE $c^{k+1} \longleftarrow \gamma c^k$.
		\STATE $n^{k+1} \longleftarrow n^k + n_{\mathrm{aug}}$.
		\STATE $k \longleftarrow k+1$.
	\ENDWHILE
	\STATE $\hat{y}\left(\mathbf{x}(\boldsymbol{\mu}) ; \boldsymbol{\theta}^*_{y}\right) \longleftarrow  \hat{y}\left(\mathbf{x}(\boldsymbol{\mu}) ; \boldsymbol{\theta}^k_{y}\right)$.
	\STATE $\hat{u}\left(\mathbf{x}(\boldsymbol{\mu}) ; \boldsymbol{\theta}^*_{u}\right) \longleftarrow  \hat{u}\left(\mathbf{x}(\boldsymbol{\mu}) ; \boldsymbol{\theta}^k_{u}\right)$.
	\ENSURE $\hat{y}\left(\mathbf{x}(\boldsymbol{\mu}) ; \boldsymbol{\theta}^*_{y}\right),  \hat{u}\left(\mathbf{x}(\boldsymbol{\mu}) ; \boldsymbol{\theta}^*_{u}\right)$.
\end{algorithmic}
\end{algorithm}

\begin{remark}
A post-processing step can be applied to continue training $\boldsymbol{\theta}^*_{y}$ until a more accurate solution of the state function (or the adjoint function) is found. That is, we can fix $\hat{u}\left(\mathbf{x}; \boldsymbol{\theta}^*_{u}\right)$ and train  $\boldsymbol{\theta}^*_{y}$ by minimizing the state loss \eqref{lossstate} . We can also fix $\boldsymbol{\theta}^*_{y}, \boldsymbol{\theta}^*_{u}$ and train $\boldsymbol{\theta}^*_{p}$ using \eqref{lossadjoint}. The initial step size $c^0$ is crucial for the convergence of Algorithm \ref{alg_AONN}. A large step size may lead to divergence of the algorithm, while a small one could result in slow convergence.
\end{remark}
\section{Comparison with other methods}
\label{sec:analysis}
Unlike solving the deterministic optimal control problems, the existence of parameters in $\mathrm{OCP}(\mb{\mu})$ causes difficulties for traditional grid-based numerical methods. A straightforward way is to convert the $\mathrm{OCP}(\mb{\mu})$ into the deterministic optimal control problem. For each realization of parameters, the $\mathrm{OCP}(\mb{\mu})$ is reduced to the following  
\begin{linenomath*}\begin{equation}\label{OCP}
    \mathrm{OCP}:\quad \left\{\begin{aligned}
    &\min _{(y,u) \in Y \times U} J(y,u),\\
    &\text { s.t. } \mathbf{F}(y,u) = 0 \ \text{ in }\Omega, \
    \text{ and } u\in U_{ad}.
    \end{aligned}\right.
\end{equation}\end{linenomath*}
The classical direct-adjoint looping (DAL) method \cite{mitter1971optimal,jameson1988aerodynamic} is a popular approach for solving this problem, where an iterative scheme is adopted to converge toward the optimal solution by solving subproblems in the KKT  system with numerical solvers (e.g. finite element methods). At each iteration in the direct-adjoint looping procedure, one first solves the governing PDE \eqref{DAL-state} and then solves the adjoint PDE \eqref{DAL-adjoint} which formulates the total gradient \eqref{DAL-grad} for the update of the control function. 
\begin{linenomath*}
\postdisplaypenalty=0
\begin{subequations}\label{DAL}
\begin{align}
&\mathbf{F}(y,u) = 0, \label{DAL-state}\\
&J_y(y,u) - \mathbf{F}^{*}_{y}(y,u)p = 0, \label{DAL-adjoint}\\
&\mathrm{d}_{u}J(y,u) = J_u(y,u) - \mathbf{F}^{*}_{u}(y,u)p. \label{DAL-grad}
\end{align}
\end{subequations}
\end{linenomath*}
Despite its effectiveness, DAL is not able to handle the $\mathrm{OCP}(\boldsymbol{\mu})$ problem, directly due to the curse of dimensionality of the discretization of $\Omega_{\mathcal{P}}$. An alternative strategy is the reduced order model (ROM) \cite{negri2013reduced}, which rely on surrogate models for parameter-dependent PDEs. The idea is that the solution of PDE for any parameter can be computed based on a few basis functions that are constructed from the solutions corresponding to some pre-selected parameters. However, it is still computationally unaffordable for ROM when the parameter-induced solution manifold does not lie on a low-dimensional subspace. 

Recently, some deep learning algorithms are used to solve the optimal control problem for a fixed parameter \cite{lu2021physics, mowlavi2022optimal}. By introducing two deep neural networks, the state function $y$ and the control function $u$ can be approximated by minimizing the following objective functional:
\begin{linenomath*}\begin{equation}
    \min_{(y,u) \in Y \times U} J(y,u) + \beta_{1} \mathbf{F}(y,u)^2 + \beta_{2} \|u-\mathbf{P}_{U_{ad}}(u)\|_{U},
    \label{eq_pinn_oc}
\end{equation}\end{linenomath*}
where two penalty terms are added, and $\beta=(\beta_{1}, \beta_{2})$ are two parameters that need tuning. As the penalty parameters increase to $+ \infty$, the solution set of this unconstrained problem approaches to the solution set of the constrained one. However,  this penalty approach has a serious drawback. On the one hand, as the penalty parameters increase, the optimal solution becomes increasingly difficult to obtain. On the other hand, the constraint is not satisfied well if the penalty parameter is small. To alleviate this difficulty, one can use hPINN \cite{lu2021physics} which employs the augmented Lagrangian method to solve \eqref{eq_pinn_oc}. However, it is still challenging to directly extend this approach to $\mathrm{OCP}(\boldsymbol{\mu})$ due to the presence of parameters. This is because it is extremely hard to optimize a series of objective functionals with a continuous range of parameters simultaneously.

\subsection{PINN for OCP(\texorpdfstring{\(\boldsymbol{\mu}\)}{μ})}
For handling the parametric optimal control problems, an extended PINN method \cite{demo2021extended} with augmented inputs is used to obtain a more accurate parametric prediction. That is, the inputs of the neural networks consist of two parts: the spatial coordinates and the parameters. The optimal solution for any parameter is approximated by a deep neural network that is obtained from solving the parameter-dependent KKT system \eqref{KKT}. In particular, when there is no restriction on the control function $u(\boldsymbol{\mu})$ (e.g., $U_{ad}(\boldsymbol{\mu})$ is the full Banach space), the KKT system is 
\begin{linenomath*}\begin{equation}\label{KKT_PINN}
    \mathcal{F}(y(\boldsymbol{\mu}),u(\boldsymbol{\mu}),p(\boldsymbol{\mu});\boldsymbol{\mu}) = \left[\begin{gathered}
     J_y(y(\boldsymbol{\mu}),u(\boldsymbol{\mu}); \boldsymbol{\mu}) - \mathbf{F}_{y}^*(y(\boldsymbol{\mu}),u(\boldsymbol{\mu}); \boldsymbol{\mu})p(\boldsymbol{\mu})\\
    \mathbf{F}(y(\boldsymbol{\mu}),u(\boldsymbol{\mu}); \boldsymbol{\mu})\\
     \mathrm{d}_u J(y(\boldsymbol{\mu}),u(\boldsymbol{\mu}); \boldsymbol{\mu})
    \end{gathered}\right]=\mb{0},
\end{equation}\end{linenomath*}
where the total gradient $d_u J$ is given in \cref{grad}. In such cases, one can use the PINN algorithm to obtain the optimal solution through minimizing the least-square loss derived from the KKT system.  Nevertheless, to apply this method to the cases where there are some additional constraints on the control function $u$, such as the box constraint \eqref{box-constraint}, one may need to introduce the Lagrange multipliers $\lambda(\boldsymbol{\mu}) = (\lambda_a(\boldsymbol{\mu}), \lambda_b(\boldsymbol{\mu}))$ corresponding to the constraints $u(\boldsymbol{\mu}) \geq u_a(\boldsymbol{\mu})$ and $u(\boldsymbol{\mu}) \leq u_b(\boldsymbol{\mu})$. For such cases, the KKT system is
\begin{linenomath*}
\postdisplaypenalty=0
\begin{subequations}\label{KKT_PINN_box}
\begin{align}
    \mathcal{F}(y(\boldsymbol{\mu}),u(\boldsymbol{\mu}),p(\boldsymbol{\mu}), \lambda(\boldsymbol{\mu}); \boldsymbol{\mu})& = \left[\begin{gathered}
     J_y(y(\boldsymbol{\mu}),u(\boldsymbol{\mu}); \boldsymbol{\mu}) - \mathbf{F}_{y}^*(y(\boldsymbol{\mu}),u(\boldsymbol{\mu}); \boldsymbol{\mu})p(\boldsymbol{\mu})\\
    \mathbf{F}(y(\boldsymbol{\mu}),u(\boldsymbol{\mu}); \boldsymbol{\mu})\\
    \mathrm{d}_u J(y(\boldsymbol{\mu}),u(\boldsymbol{\mu}); \boldsymbol{\mu}) - \lambda_a(\boldsymbol{\mu}) + \lambda_b(\boldsymbol{\mu})\\
    \lambda_a(\boldsymbol{\mu})(u_a(\boldsymbol{\mu})-u(\boldsymbol{\mu}))\\
    \lambda_b(\boldsymbol{\mu})(u_b(\boldsymbol{\mu})-u(\boldsymbol{\mu}))
    \end{gathered}\right]=\mb{0}, \label{KKT_PINN_box1}\\ 
    \text{ and } \quad  &
    \left\{\begin{aligned}
    u_a(\boldsymbol{\mu})\leq u(\boldsymbol{\mu}) \leq u_b(\boldsymbol{\mu}),\\
    \lambda_a(\boldsymbol{\mu}) \geq0, \lambda_b(\boldsymbol{\mu}) \geq 0.
    \end{aligned}\right. \label{KKT_PINN_box2}
\end{align}
\end{subequations}
\end{linenomath*}
Applying the framework of PINN to solve the system \eqref{KKT_PINN_box} needs to deal with several penalty terms in the loss function including the penalties of equality terms \eqref{KKT_PINN_box1} and inequality terms \eqref{KKT_PINN_box2}, leading to an inaccurate solution even for the problem with fixed parameters, which will be presented in the next section. In addition, the extra constraint of $U_{ad}(\boldsymbol{\mu})$ often introduces inequality terms and nonlinear terms and brings more singularity to the optimal control function \cite{ali2020reduced, bader2016certified}, which limits the application of PINN for solving $\mathrm{OCP}(\boldsymbol{\mu})$ with control constraints.

\subsection{PINN+Projection for OCP(\texorpdfstring{\(\boldsymbol{\mu}\)}{μ})}

To find a better baseline for comparison, we propose to improve the performance of PINN by introducing a projection operator.
In this way, the KKT system \eqref{KKT} can be reformulated to a more compactly stated condition \cite{nocedal1999numerical}:
\begin{linenomath*}\begin{equation}\label{projectedKKT}
    \mathcal{F}(y(\boldsymbol{\mu}),u(\boldsymbol{\mu}),p(\boldsymbol{\mu}),c; \boldsymbol{\mu}) = \left[\begin{gathered}
    J_y(y(\boldsymbol{\mu}),u(\boldsymbol{\mu}); \boldsymbol{\mu}) - \mathbf{F}_{y}^*(y(\boldsymbol{\mu}),u(\boldsymbol{\mu}); \boldsymbol{\mu})p(\boldsymbol{\mu})\\
    \mathbf{F}(y(\boldsymbol{\mu}),u(\boldsymbol{\mu}); \boldsymbol{\mu})\\
     u(\boldsymbol{\mu}) - \mathbf{P}_{U_{ad}(\boldsymbol{\mu})}\left(u(\boldsymbol{\mu})-c \mathrm{d}_{u}J(y(\boldsymbol{\mu}),u(\boldsymbol{\mu}); \boldsymbol{\mu})\right)
    \end{gathered}\right]=\mb{0},
\end{equation}\end{linenomath*}
where $c$ could be any positive number. Note that choosing an appropriate $c$ can accelerate the convergence of the algorithm. For example, the classic way is to choose $c=1/\alpha$ for canceling out the control function $u$ inside the projection operator, where $\alpha$ is the coefficient of the Tikhonov regularization term (see the experiment in Section~\ref{sec_test1}).
The complementary conditions and inequalities caused by the control constraints are avoided in \cref{projectedKKT}, thus significantly reducing the difficulty of optimization.
In this paper, we call the method PINN+Projeciton, which combines the projection strategy with the KKT system to formulate the PINN residual loss.
Although PINN+Projection alleviates the solving difficulty brought by control constraints to PINN, it still has limitations on nonsmooth optimal control problems.
For nonsmooth optimization such as sparse $L_1$-minimization, the KKT system can no longer be described by \eqref{projectedKKT} because of the nondifferentiable property of the $L_1$-norm \cite{de2015numerical}. 
Instead, the dual multiplier for the $L_1$-cost term is required. The difficulty arises from the third nonsmooth variational equation of \eqref{projectedKKT}, which makes the neural network difficult to train. 
AONN reduces this difficulty by leveraging the update scheme in the DAL method without the implicit variational equation in \eqref{projectedKKT}. Numerical results also show that the KKT system \eqref{projectedKKT} cannot be directly used to formulate the loss functions of neural networks. 
Such results of $L_1$-minimization involved in $\mathrm{OCP}$($\mb{\mu}$) (see Section~\ref{sec_test5_sparse}) strongly suggest that AONN is a more reliable and efficient framework. 

The proposed AONN method has all the advantages of the aforementioned approaches while avoiding their drawbacks. By inheriting the structure of DAL, the AONN method can obtain an accurate solution through solving the KKT system in an alternative minimization iterative manner. So it does not require the Lagrange multipliers corresponding to the additional control constraints and thus can reduce the storage cost as well as improve the accuracy. Moreover,  AONN can accurately approximate the optimal solutions of parametric optimal control problems for any parameter and can be generalized to cases with high-dimensional parameters.

\section{Numerical study}
\label{sec:results}

In this section, we present results of five numerical experiments to illustrate the effectiveness of AONN, where different types of PDE constraints, objective functionals and control constraints under different parametric settings are studied. In the following, AONN is first validated by solving OCP, and further applied to solving OCP($\mb{\mu})$ with continuous parameters changing over a specific interval. For comparison purposes, we also use the PINN method and the PINN+Projecton method to solve OCP($\mb{\mu}$).
We employ the ResNet model \cite{he2016deep} with sinusoid activation functions to build the neural networks for AONN and other neural network based algorithms. Unless otherwise specified, the quasi Monte-Carlo method is used to generate collocation points from $\Omega_{\mathcal{P}}$ by calling the SciPy module \cite{virtanen2020scipy}. Analytical length factor functions (see \eqref{eq-pfnn}) are constructed for all test problems to make the approximate solution naturally satisfy Dirichlet boundary conditions. The training of neural networks is performed on a Geforce RTX 2080 GPU with PyTorch 1.8.1. The Broyden–Fletcher–Goldfarb–Shanno (BFGS) algorithm with a strong Wolfe line search strategy is used to update the neural network parameters to speed up the convergence, where the maximal number of iterations for BFGS is set to $100$. 

\subsection{Test 1: Optimal control for the semilinear elliptic equations}\label{sec_test1}
We start with the following nonparametric optimal control problem:
\begin{linenomath*}\begin{equation}\label{test1}
    \left\{\begin{aligned} 
    &\min_{y,u}  J(y, u):=\frac{1}{2}\left\|y-y_{d}\right\|_{L_{2}(\Omega)}^{2}+\frac{\alpha}{2}\|u\|_{L_{2}(\Omega)}^{2} ,\\
    &\text { subject to }  
    \left\{\begin{aligned}
    -\Delta y + y^3&=u+f &&\text{ in } \Omega \\ 
    y &=0  &&\text { on } \partial \Omega,
    \end{aligned}\right.\\ 
    &\text{and}\quad u_a \leq u \leq u_b  \quad \text { a.e. in } \Omega.\\
    \end{aligned}\right.
\end{equation}\end{linenomath*}
The total derivative of $J$ with respect to $u$ is $\mathrm{d}_{u}J(y,u) = \alpha u + p$, where $p$ is the solution of the corresponding adjoint equation:
\begin{linenomath*}\begin{equation}
\label{test1-adjoint}
\left\{\begin{aligned}
-\Delta p + 3py^2 &= y-y_d &&\text{ in } \Omega, \\ 
p &=0  &&\text { on } \partial \Omega.
\end{aligned}\right.
\end{equation}\end{linenomath*}
We take the same configuration as in ref.\cite{gong2015multilevel}, where $\Omega=(0,1)^2$, $\alpha=0.01$, $u_a=0$, and $u_b=3$. The analytical optimal solution is given by
\begin{linenomath*}\begin{equation}\label{case1}
    \begin{aligned}
    &y^{*} = \sin \left(\pi x_{1}\right) \sin \left(\pi x_{2}\right),\\
    &u^{*} = \mathbf{P}_{[u_a,u_b]}(2\pi^{2}y^{*}),\\
    &p^{*} = -2\alpha\pi^{2}y^{*},\\
    \end{aligned}
\end{equation}\end{linenomath*}
where $\mathbf{P}_{[u_a,u_b]}$ is the pointwise projection operator onto the interval $[u_a,u_b]$. The desired state $y_d = (1+4\pi^4\alpha)y^{*}-3y^{*2}p^{*}$ and the source term $f = 2\pi^2 y^{*} + y^{*3} - u^{*}$ are given to satisfy the state equation and the adjoint equation.

\begin{figure}[!htbp]
  \centering
  \subfigure[]
  {\includegraphics[width=0.38\textwidth]{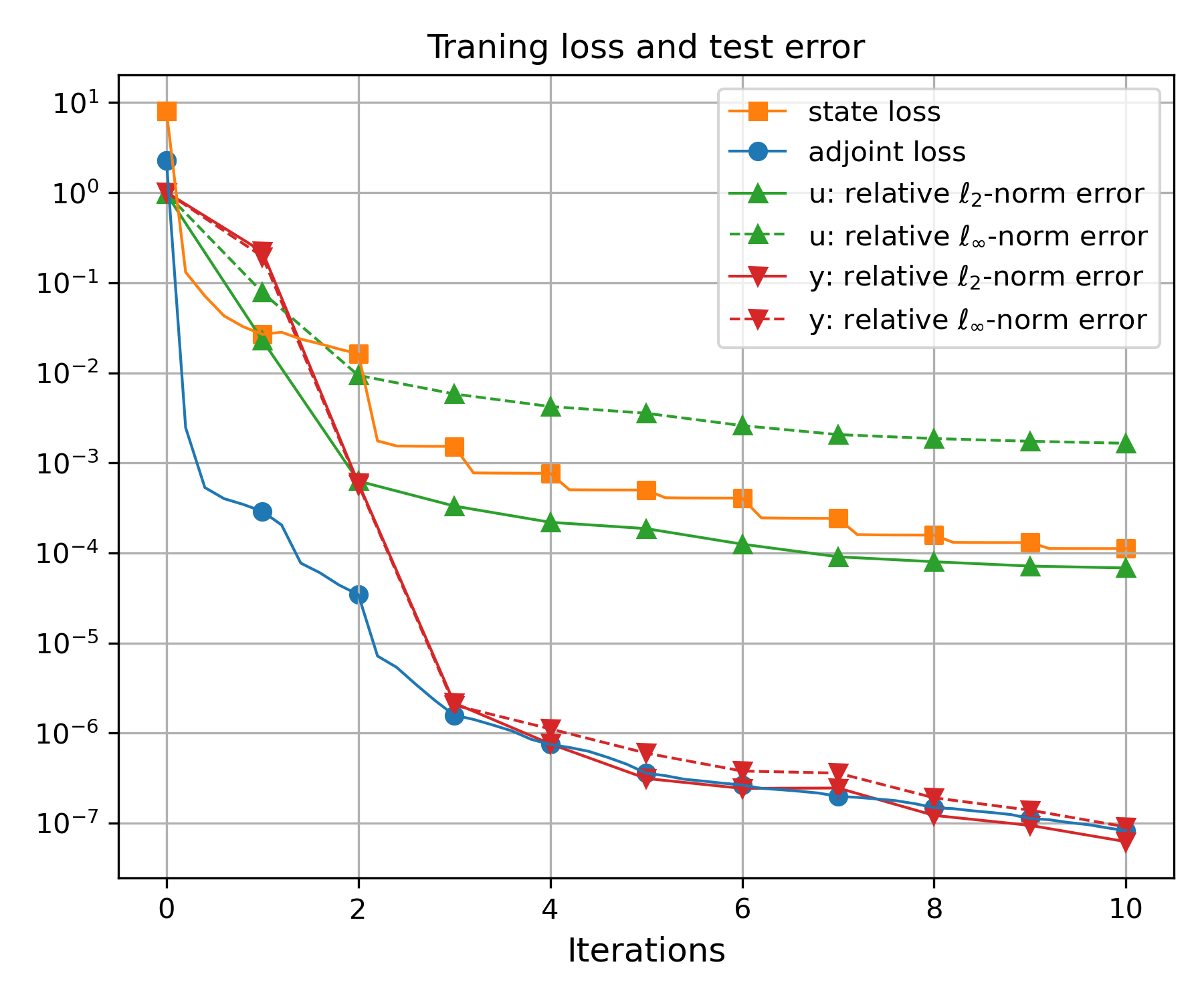}\label{fig:test1-ocp-1}}
  \subfigure[]
  {\includegraphics[width=0.58\textwidth]{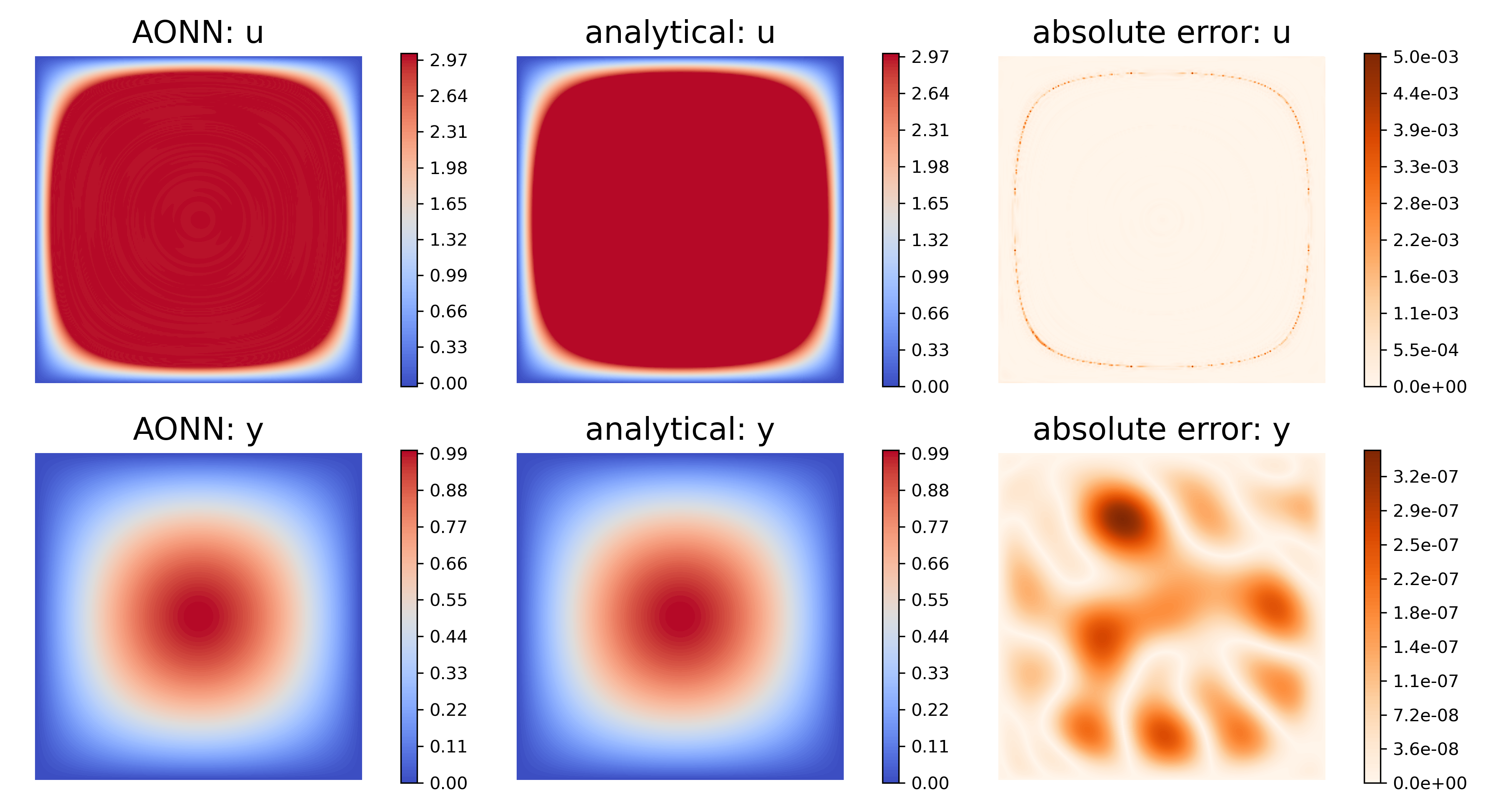}\label{fig:test1-ocp-2}}
  \caption{Test 1: training loss and test error of problem \ref{test1}. Test error is evaluated at 256$\times$256 grid points. (a) Loss behaviour measured in terms of \eqref{lossstate}-\eqref{lossadjoint}, and  test errors in both $\ell_2$-norm and $\ell_{\infty}$-norm during training process. (b) The AONN solution and its absolute errors compared with the analytical solution.}
\end{figure}

To solve the optimal control problem with AONN, we construct three networks $\hat{y}_{I}\left(\mathbf{x}(\boldsymbol{\mu}); \boldsymbol{\theta}_{y_I}\right), \hat{p}_{I}\left(\mathbf{x}(\boldsymbol{\mu}); \boldsymbol{\theta}_{p_I}\right)$ and $\hat{u}\left(\mathbf{x}(\boldsymbol{\mu}); \boldsymbol{\theta}_{u}\right)$, whose network structures are all comprised of two ResNet blocks, each of which contains two fully connected layers with $15$ neurons and a residual connection. We randomly sample $N = 4096$ points inside $\Omega$ to form the training set.
A uniform meshgrid with size $256\times256$ in $\Omega$ is generated for testing and visualization. We use fixed step size and training epochs in subproblems, i.e. $c^k \equiv 1/\alpha = 100, n^k\equiv500$. The loss behavior and the relative error $\|u - u^*\|/\|u^*\|$ with $\ell_2$-norm and $\ell_{\infty}$-norm are reported in Figure~\ref{fig:test1-ocp-1}, while Figure~\ref{fig:test1-ocp-2} evaluates the difference between the AONN solution and the analytical solution.
As reported in ref.\cite{gong2015multilevel}, to achieve the error $1 \times 10^{-4}$ in the $\ell_2$ sense requires $7733$ degrees of freedom with the finite element method, while the AONN method needs only $781$ neural network parameters to approximate the control function.

\subsection{Test 2: Optimal control for the semilinear elliptic equations with control constraint parametrization}\label{sec_test2}
We then consider the same optimal control problem with control constraint parametrization. 
The control constraint upper bound $u_b$ is set to be a continuous variable $\boldsymbol{\mu}$ ranging from $3$ to $20$ instead of a fixed number. Thus \eqref{test1} actually constructs a series of optimal control problems and the optimal solutions \eqref{case1} are dependent on $\boldsymbol{\mu}$. We now verify whether the all-at-once solutions can be obtained by AONN when $\boldsymbol{\mu}$ changes continuously over the interval $3\leq \boldsymbol{\mu} \leq 20$. We seek optimal $y(\boldsymbol{\mu}), u(\boldsymbol{\mu})$ defined by the following problem:
\begin{linenomath*}\begin{equation}
    \left\{\begin{aligned} 
    &\min_{y(\boldsymbol{\mu}),u(\boldsymbol{\mu})}  J(y(\boldsymbol{\mu}), u(\boldsymbol{\mu})):=\frac{1}{2}\left\|y(\boldsymbol{\mu})-y_{d}(\boldsymbol{\mu})\right\|_{L_{2}(\Omega)}^{2}+\frac{\alpha}{2}\|u(\boldsymbol{\mu})\|_{L_{2}(\Omega)}^{2}, \\
    &\text { subject to }  
    \left\{\begin{aligned}
    -\Delta y(\boldsymbol{\mu}) + y(\boldsymbol{\mu})^3&=u(\boldsymbol{\mu})+f(\boldsymbol{\mu}) &&\text{ in } \Omega \\ 
    y(\boldsymbol{\mu}) &=0  &&\text { on } \partial \Omega,
    \end{aligned}\right.\\ 
    &\text{and}\quad u_a \leq u(\boldsymbol{\mu}) \leq \boldsymbol{\mu} \quad \text { a.e. in } \Omega.\\
    \end{aligned}\right.
    \label{test1-pocp}
\end{equation}\end{linenomath*}
To naturally satisfy the homogeneous Dirichlet boundary conditions in the state equation and the adjoint equation, three neural networks for approximating the AONN solutions of $\mathrm{OCP}(\boldsymbol{\mu})$  \eqref{test1-pocp} are defined as follows:
\begin{linenomath*}\begin{equation}
    \begin{aligned}
    \hat{y}\left(\mathbf{x}(\boldsymbol{\mu}) ; \boldsymbol{\theta}_{y_I}\right)&=\ell(\mathbf{x}) \hat{y}_{I}\left(\mathbf{x}(\boldsymbol{\mu}) ; \boldsymbol{\theta}_{y_I}\right),\\
    \hat{p}\left(\mathbf{x}(\boldsymbol{\mu}) ; \boldsymbol{\theta}_{p_I}\right)&=\ell(\mathbf{x}) \hat{p}_{I}\left(\mathbf{x}(\boldsymbol{\mu}) ; \boldsymbol{\theta}_{p_I}\right),\\
    \hat{u}\left(\mathbf{x}(\boldsymbol{\mu}) ; \boldsymbol{\theta}_{u}\right)&= \hat{u}_{I}\left(\mathbf{x}(\boldsymbol{\mu}) ; \boldsymbol{\theta}_{u}\right),
    \end{aligned}\\
\end{equation}\end{linenomath*}
where the length factor function is formed by
\begin{linenomath*}\begin{equation}
    \ell(\mathbf{x}) = x_0(1 - x_0)x_1(1 - x_1).
\end{equation}\end{linenomath*}

The network structures of $\hat{y}_{I}\left(\mathbf{x}(\boldsymbol{\mu}); \boldsymbol{\theta}_{y_I}\right), \hat{p}_{I}\left(\mathbf{x}(\boldsymbol{\mu}); \boldsymbol{\theta}_{p_I}\right)$ and $\hat{u}\left(\mathbf{x}(\boldsymbol{\mu}); \boldsymbol{\theta}_{u}\right)$ are the same as those of the previous test except for the input dimension being $3$ and the number of neurons in each hidden layer being $20$, resulting in $1361$ undecided parameters. To evaluate the loss, we sample $N = 20480$ points in the spatio-parametric space $\Omega_{\mathcal{P}}$. We keep the same step size $c^k$ and training epochs $n^k$ as the previous test, and perform $N_{\mathrm{iter}}=20$ iterations  until Algorithm \ref{alg_AONN} converges. The test errors are computed on the uniform meshgrid with size $256 \times256$ for each realization of $\boldsymbol{\mu}$.

\begin{figure}[htbp]
\centering
\includegraphics[width=14cm]{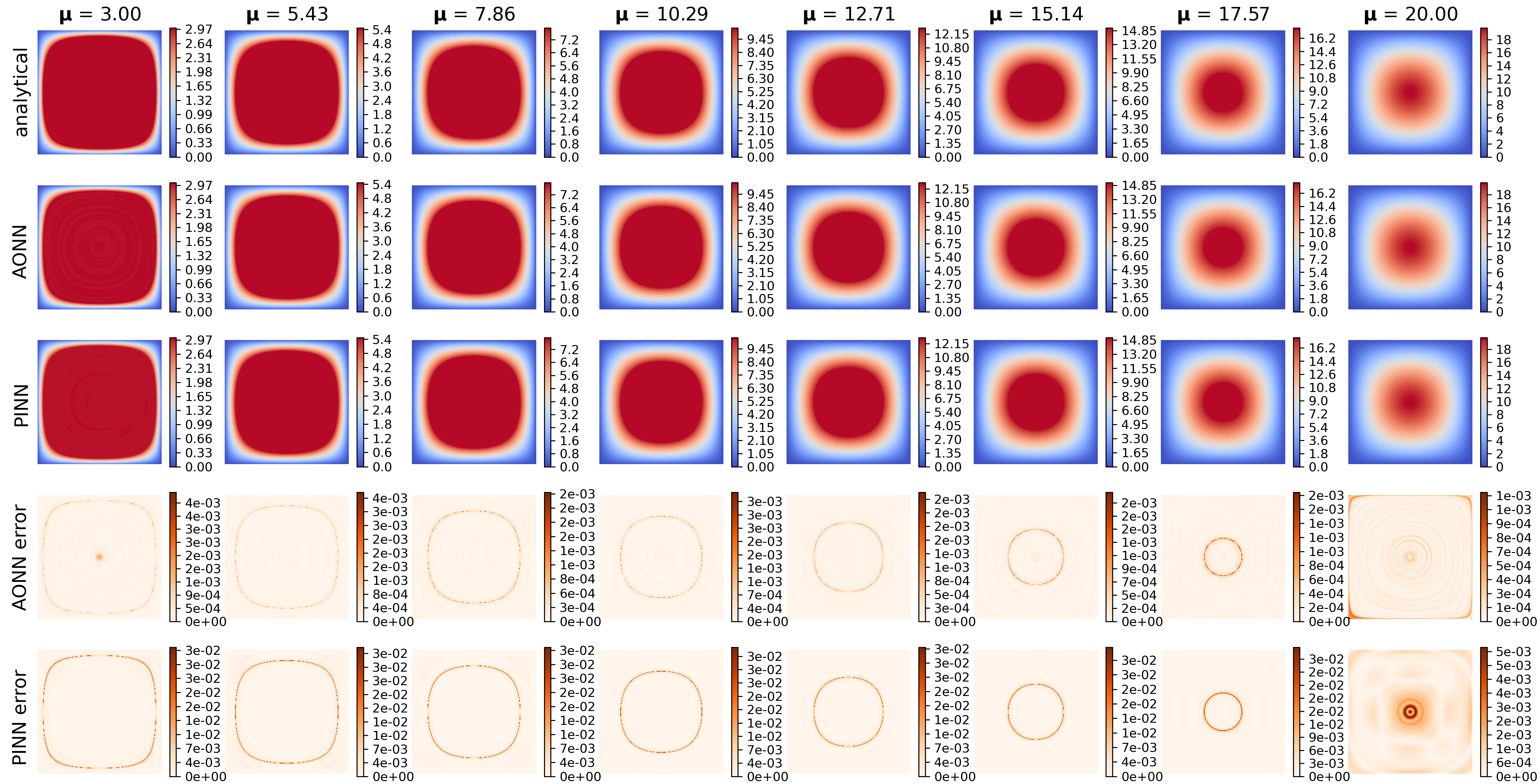}
\caption{Test 2: the control solutions $u(\boldsymbol{\mu})$ of AONN and PINN with eight realizations of $\mb{\mu}\in[3,20]$, and their absolute errors.}
\label{fig:test1-pocp}
\end{figure}

In Figure \ref{fig:test1-pocp}, we plot the analytical solutions, the AONN solutions and the PINN solutions for eight equidistant realizations of $\mb{\mu}$, where it can be seen that the AONN solutions are better than the PINN solutions in the sense of absolute error.
Looking more closely, the large errors are concentrated around the location of singularity of $u$, i.e., the curve of active constraints $\{\mathbf{x} : u(\mathbf{x}(\boldsymbol{\mu}))=\boldsymbol{\mu}\}$, except for the case $\boldsymbol{\mu}=20$ where the inequality constraint is nonactive, keeping the smoothness of the optimal control function. Note that adaptive sampling strategies~\cite{tang2021das, adda_2022, gao2022failure} may be used to improve the accuracy in the singularity region, which will be left for future study.

\subsection{Test 3: Optimal control for the Navier-Stokes equations with physical parametrization}\label{sec_test3}
The next test case is the  parametric optimal control problem
\begin{linenomath*}\begin{equation}
   \min _{y(\boldsymbol{\mu}), u(\boldsymbol{\mu})} J(y(\boldsymbol{\mu}),u(\boldsymbol{\mu}))=\frac{1}{2}\left\|y(\boldsymbol{\mu})-y_{d}(\boldsymbol{\mu})\right\|_{L_2(\Omega)}^{2}+\frac{1}{2}\|u(\boldsymbol{\mu})\|_{L_2(\Omega)}^{2},
\end{equation}\end{linenomath*}
subject to the following steady-state incompressible Navier-Stokes (NS) equations:
\begin{linenomath*}\begin{equation}\label{NS}
\left\{
\begin{aligned}
-\boldsymbol{\mu} \Delta y(\boldsymbol{\mu})+(y(\boldsymbol{\mu}) \cdot \nabla) y(\boldsymbol{\mu})+\nabla p(\boldsymbol{\mu}) &=u(\boldsymbol{\mu}) + f(\boldsymbol{\mu}) & & \text { in } \Omega, \\
\operatorname{div} y(\boldsymbol{\mu}) &=0 & & \text { in } \Omega, \\
y(\boldsymbol{\mu}) &=0 & & \text { on } \partial \Omega,
\end{aligned}
\right.
\end{equation}\end{linenomath*}
in $\Omega=(0,1)^2$ with parameter $\boldsymbol{\mu}$ representing the reciprocal of the Reynolds number.
Note that the nonparametric problems without control constraint for $\boldsymbol{\mu} = 0.1$ and $\boldsymbol{\mu} = 1.0$ were studied in refs.\cite{lowe2011projection,wachsmuth2006optimal}. We set the physical parameter $\boldsymbol{\mu}\in[0.1, 100]$ and in addition, we consider the following constraint for $u(\boldsymbol{\mu})=(u_1(\boldsymbol{\mu}), u_2(\boldsymbol{\mu}))$:
\begin{linenomath*}\begin{equation}\label{u-constrain}
    u_1(\boldsymbol{\mu})^2 + u_2(\boldsymbol{\mu})^2 \leq r^2,
\end{equation}\end{linenomath*}
with $r = 0.2$, posing additional challenges to this problem.
The desired state $y_{d}(\boldsymbol{\mu})$ and the source term $f(\boldsymbol{\mu})$ are given in advance to ensure that the analytical solution of the above OCP($\mb{\mu}$) is given by 
\begin{linenomath*}
\begin{equation*}
\begin{aligned}
y^{*}(\boldsymbol{\mu})&=e^{-0.05 \boldsymbol{\mu}}\left(\begin{array}{c}
\sin ^{2} \pi x_{1} \sin \pi x_{2} \cos \pi x_{2} \\
-\sin ^{2} \pi x_{2} \sin \pi x_{1} \cos \pi x_{1}
\end{array}\right),\\
\lambda^{*}(\boldsymbol{\mu})&=\left(e^{-0.05 \boldsymbol{\mu}}-e^{-\boldsymbol{\mu}}\right)\left(\begin{array}{c}
\sin ^{2} \pi x_{1} \sin \pi x_{2} \cos \pi x_{2} \\
-\sin ^{2} \pi x_{2} \sin \pi x_{1} \cos \pi x_{1},
\end{array}\right).
\end{aligned}
\end{equation*}
\end{linenomath*}
The adjoint equation is specified as
\begin{linenomath*}\begin{equation}\label{NSadjoint}
\left\{
\begin{aligned}
-\boldsymbol{\mu} \Delta \lambda(\boldsymbol{\mu})-(y(\boldsymbol{\mu}) \cdot \nabla) \lambda(\boldsymbol{\mu})+(\nabla y(\boldsymbol{\mu}))^{T} \lambda(\boldsymbol{\mu})+\nabla \nu(\boldsymbol{\mu}) &=y(\boldsymbol{\mu})-y_{d}(\boldsymbol{\mu}) & & \text { in } \Omega, \\
\operatorname{div} \lambda(\boldsymbol{\mu}) &=0 & & \text { in } \Omega, \\
\lambda(\boldsymbol{\mu}) &=0 & & \text { on } \partial \Omega,
\end{aligned} 
\right.
\end{equation}\end{linenomath*}
where $\lambda(\boldsymbol{\mu})$ denotes the adjoint velocity and  $\nu(\boldsymbol{\mu})$ denotes the adjoint pressure. The optimal pressure and adjoint pressure $p^{*}(\boldsymbol{\mu}), \nu^{*}(\boldsymbol{\mu})$ are both zero. In order to satisfy the state equation and the adjoint equation,  $y_{d}(\boldsymbol{\mu})$ and $f(\boldsymbol{\mu})$ are chosen as
\begin{linenomath*}\begin{equation}
\begin{aligned}
f(\boldsymbol{\mu}) &= -\boldsymbol{\mu} \Delta y^{*}(\boldsymbol{\mu})+(y^{*}(\boldsymbol{\mu}) \cdot \nabla) y^{*}(\boldsymbol{\mu})-u^{*}(\boldsymbol{\mu}), \\
y_{d}(\boldsymbol{\mu})&=y^{*}(\boldsymbol{\mu})-\left(-\boldsymbol{\mu} \Delta \lambda^{*}(\boldsymbol{\mu})+(y^{*}(\boldsymbol{\mu}) \cdot \nabla) \lambda^{*}(\boldsymbol{\mu})-(\nabla y^{*}(\boldsymbol{\mu}))^{T} \lambda^{*}(\boldsymbol{\mu})\right).
\end{aligned}
\end{equation}\end{linenomath*}
It is easy to check that the optimal control is $u^{*}(\boldsymbol{\mu})=\mathbf{P}_{B(0,r)}(\lambda^{*}(\boldsymbol{\mu}))$, where $B(0,r)$ is a ball centered at the origin of radius $r$. The state equation \eqref{NS} and the adjoint equation \eqref{NSadjoint} together with the variational inequality where $\mathrm{d}_{u}J(\boldsymbol{\mu})= u(\boldsymbol{\mu}) - \lambda(\boldsymbol{\mu})$ formulate the optimality system.

For AONN, we use a neural network to approximate $y$,  and it is constructed by two ResNet blocks, each of which contains two fully connected layers with $20$ units and a residual connection,  resulting in $1382$ parameters. The neural network for approximating $p$ has two ResNet blocks built by two fully connected layers with $10$ units, resulting in $381$ parameters. The architectures of the neural networks for $\lambda$ and $\nu$ are the same as those of $y$ and $p$ respectively. 
We select $N=20000$ randomly sampled points in the spatio-parametric space $\Omega_{\mathcal{P}}$. The maximum iteration number in Algorithm \ref{alg_AONN} is set to $N_{\mathrm{iter}} = 300$ and the step size is $c^k\equiv c^0=1.0$. 
We choose an initial training epoch $n^0 = 200$ and increase it by $n_{\mathrm{aug}}=100$ after every $100$ iterations.
For the PINN method, the architectures of the neural networks are the same as those of AONN except for adding another neural network for $\zeta(\boldsymbol{\mu})$ to satisfy the following KKT system:
\begin{linenomath*}\begin{equation}\label{test3-KKT}
    \left\{\begin{array}{l}
    \text{state equation}\, \eqref{NS},\\
    \text{adjoint equation}\, \eqref{NSadjoint},\\
    u_1(\boldsymbol{\mu}) - \lambda_1(\boldsymbol{\mu}) + 2u_1(\boldsymbol{\mu}) \zeta(\boldsymbol{\mu}) = 0,\\
    u_2(\boldsymbol{\mu}) - \lambda_2(\boldsymbol{\mu}) + 2u_2(\boldsymbol{\mu}) \zeta(\boldsymbol{\mu}) = 0,\\
    (u_1(\boldsymbol{\mu})^2 + u_2(\boldsymbol{\mu})^2 - r^2)\zeta(\boldsymbol{\mu}) = 0,\\
    u_1(\boldsymbol{\mu})^2 + u_2(\boldsymbol{\mu})^2 \leq r^2, \zeta(\boldsymbol{\mu}) \geq 0,\\
    \end{array}\right.\\
\end{equation}\end{linenomath*}
 where $\zeta(\boldsymbol{\mu})$ is the Lagrange multiplier of the control constraint \eqref{u-constrain}.

\begin{figure}[htbp]
\centering 
\includegraphics[height=5.5 cm]{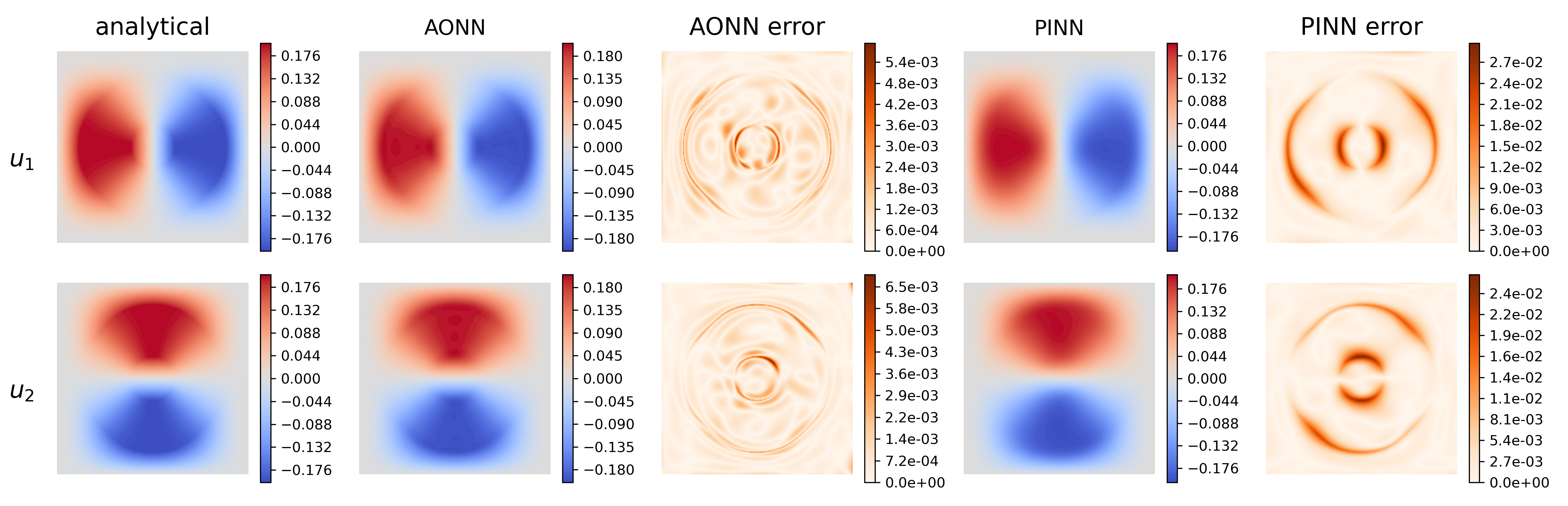}
\caption{Test 3: optimal solutions of the state function $y=(y_1,y_2)$ and the control function $u=(u_1,u_2)$ obtained by AONN and PINN, and their absolute errors for a given parameter $\boldsymbol{\mu}=10$.}
\label{fig:test2-compare}
\end{figure}

\begin{figure}[htbp]
\centering
\includegraphics[height=5cm]{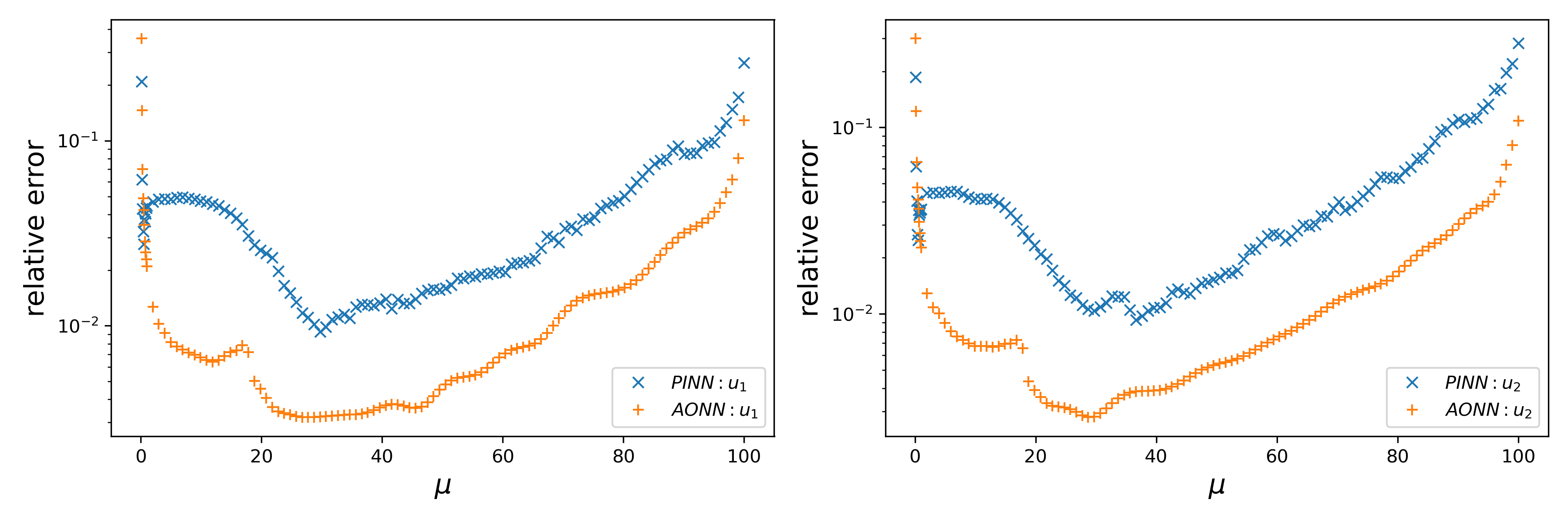}
\caption{Test 3: the relative errors (in the $\ell_2$-norm sense) of AONN and PINN for the two components of $u(\boldsymbol{\mu})=(u_1(\boldsymbol{\mu}), u_2(\boldsymbol{\mu}))$. The relative errors are computed on the $256\times 256$ meshgrid for each fixed parameter $\boldsymbol{\mu}$.}
\label{fig:test2-mu}
\end{figure}

We compare the solutions of AONN with those obtained using PINN and plot their absolute errors in Figure~\ref{fig:test2-compare}, where it shows the control function $u=(u_1,u_2)$ for a representative parameter $\boldsymbol{\mu}=10$. 
From the figure, it can be seen that AONN can obtain a more accurate optimal control function than that of PINN, even when the training of PINN costs more epochs than that of AONN. Also, the quadratic constraint is not satisfied well for the PINN solution because there are more penalties from the KKT system \cref{test3-KKT} for the PINN loss.
We compute the relative error $\|u - u^*\|/\|u^*\|$ on a uniform $256\times 256$ meshgrid for each parameter $\boldsymbol{\mu}$ and plot the results in Figure~\ref{fig:test2-mu}. For most of the parameters, the relative errors of the AONN solutions are smaller than that of PINN, indicating that AONN is more effective and efficient than PINN in solving parametric optimal control problems.
Note that this problem becomes harder when the parameter $\boldsymbol{\mu}$ gets smaller \cite{wachsmuth2006optimal}. In particular, the relative errors of AONN and PINN are both large as $\boldsymbol{\mu}$ closes to $0.1$.

\subsection{Test 4: Optimal control for the Laplace equation with geometrical parametrization}\label{sec_test4_geo}
In this test case, we are going to solve the following parametric optimal control problem: 
\begin{linenomath*}\begin{equation}
    \label{eq:test3}
    \left\{\begin{aligned} 
    &\min _{y(\boldsymbol{\mu}), u(\boldsymbol{\mu})} J\left(y(\boldsymbol{\mu}), u(\boldsymbol{\mu})\right)=\frac{1}{2}\left\|y(\boldsymbol{\mu})-y_{d}(\boldsymbol{\mu})\right\|_{L_{2}\left(\Omega(\boldsymbol{\mu})\right)}^{2}+\frac{\alpha}{2}\left\|u(\boldsymbol{\mu})\right\|_{L_{2}\left(\Omega(\boldsymbol{\mu})\right)}^{2}, \\
    &\text { subject to }  
    \left\{\begin{aligned}
    -\Delta y(\boldsymbol{\mu})&=u(\boldsymbol{\mu}) &&\text { in } \Omega(\boldsymbol{\mu}), \\ 
    y(\boldsymbol{\mu})&=1 &&\text { on }\partial \Omega(\boldsymbol{\mu}),\\
    \end{aligned}\right.\\ 
    &\text{and} \quad u_a \leq u(\boldsymbol{\mu}) \leq u_b \quad \text { a.e. in } \Omega(\boldsymbol{\mu}),
    \end{aligned}\right.
\end{equation}\end{linenomath*}
where $\boldsymbol{\mu}=(\mu_1, \mu_2)$ represents the geometrical and desired state parameters. The parametric computational domain is $\Omega(\boldsymbol{\mu}) = ([0,2]\times[0,1])\backslash B((1.5,0.5), \mu_1)$ and the desired state is given by
\begin{linenomath*}\begin{equation}
y_d(\boldsymbol{\mu}) = 
\begin{cases}
1 & \text { in } \Omega_{1}=[0,1]\times[0,1], \\ 
\mu_2 & \text { in }\Omega_{2}(\boldsymbol{\mu})=([1,2]\times[0,1])\backslash B((1.5,0.5),\mu_1),
\end{cases}
\end{equation}\end{linenomath*}
where $B((1.5,0.5), \mu_1)$ is a ball of radius $\mu_1$ with center $(1.5,0.5)$. We set $\alpha=0.001$ and the parameter interval to be $\boldsymbol{\mu}\in \mathcal{P}=[0.05,0.45]\times[0.5,2.5]$.

\begin{figure}[!htbp]
  \centering
  \subfigure[]
  {\includegraphics[height=0.29\textwidth]{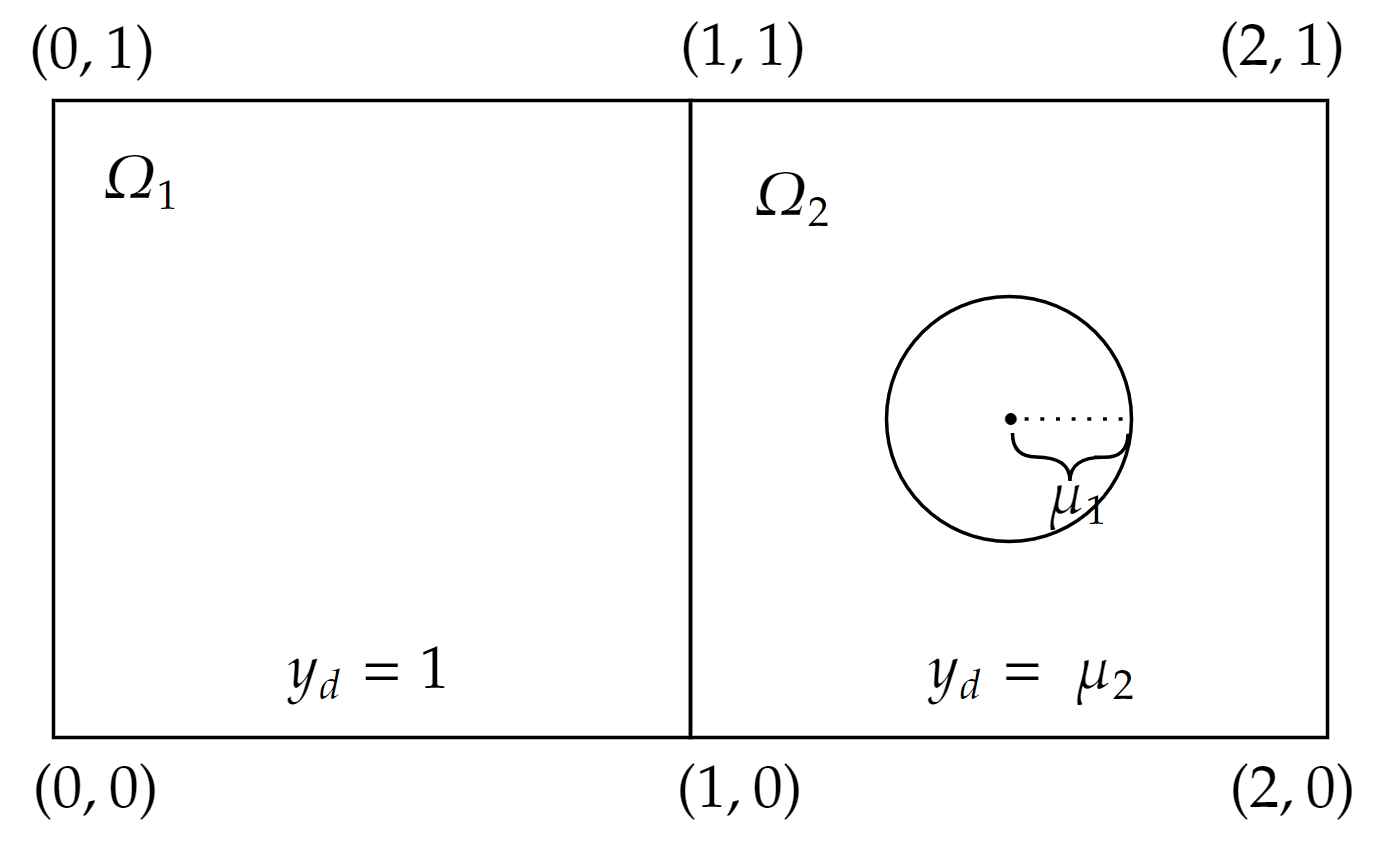}  \label{fig:test3-data1}}
  \subfigure[]
  {\includegraphics[height=0.3\textwidth]{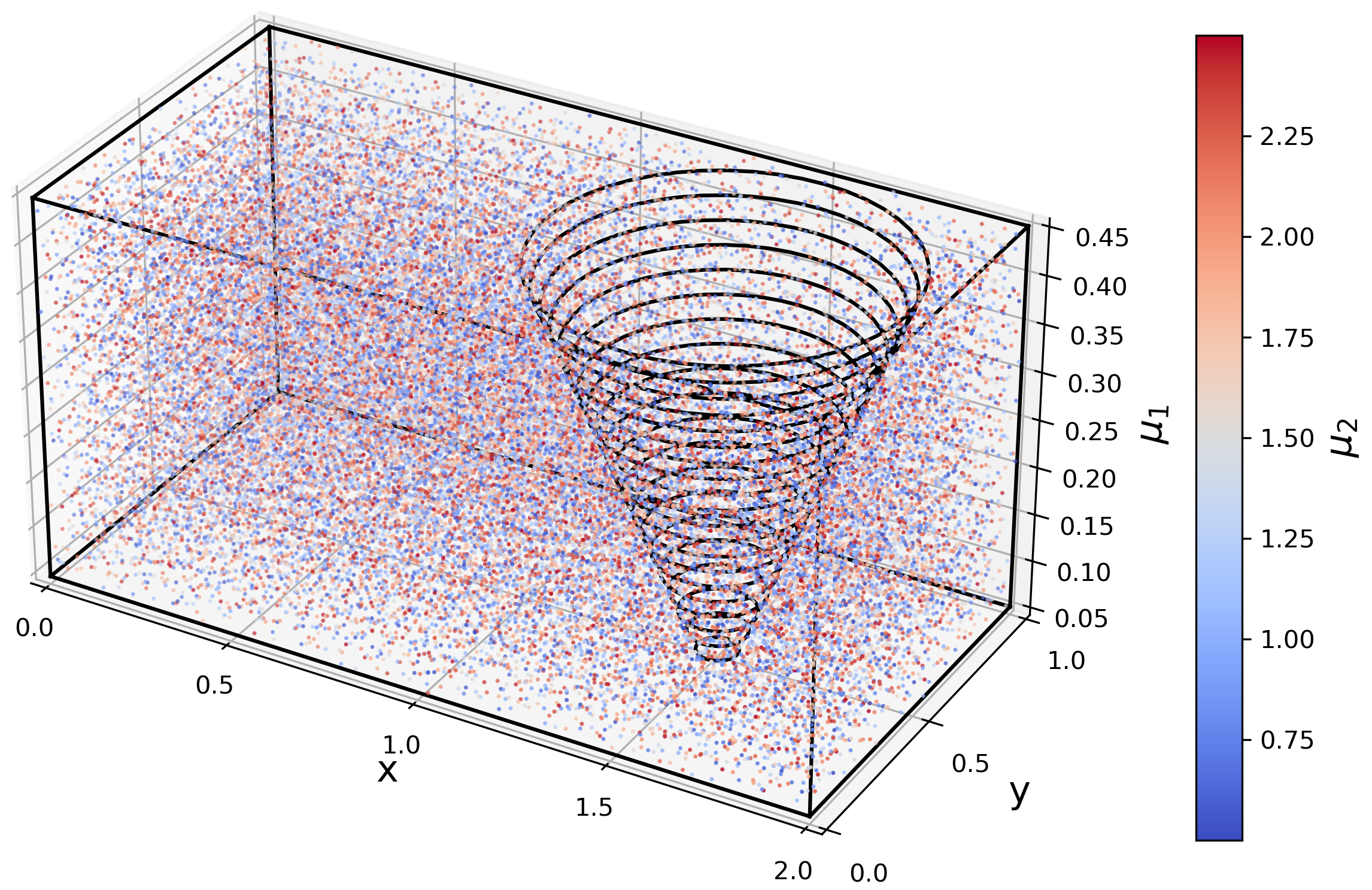}  \label{fig:test3-data2}}
  \caption{Test 4: (a) The parametric computational domain $\Omega(\boldsymbol{\mu})$. (b) $N = 40000$ training collocation points sampled in  $\Omega_{\mathcal{P}}$ (there are no points inside the frustum).}
  \label{fig:test3-data}
\end{figure}

This test case is inspired by the literature \cite{negri2013reduced, karcher2018certified} that involve the application of local hyperthermia treatment of cancer. In such case, it is expected to achieve a certain temperature field in the tumor area and another temperature field in the non-lesion area through the heat source control. The circle cut out from the rectangular area represents a certain body organ, and by using AONN we aim to obtain all-at-once solutions of the optimal heat source control for different expected temperature fields and different organ shapes. 
In particular, we consider a two-dimensional model problem corresponding to the hyperthermia cancer treatment.
One difficulty of this problem is the geometrical parameter $\mu_1$ that leads to various computational domains, which causes difficulties in applying traditional mesh-based numerical methods. In the AONN framework, we can solve this problem by sampling in the spatio-parametric space:
\begin{linenomath*}
\postdisplaypenalty=0
\begin{align*}
    \Omega_{\mathcal{P}}=\{(x_0,x_1,\mu_1,\mu_2) |& 0\leq x_0\leq2,\, 0\leq x_1\leq1,\, 0.05\leq \mu_1\leq0.45,\, 0.5\leq \mu_2\leq2.5,\,
    (x_0-1.5)^2+(x_1-0.5)^2\geq \mu_1^2\}.
\end{align*}
\end{linenomath*}
The computational domain $\Omega(\boldsymbol{\mu})$ as well as the $40000$ training points are given in Figure~\ref{fig:test3-data1} and Figure~\ref{fig:test3-data2}. 

The state neural network $\hat{y}$ is constructed by $ \hat{y}\left(\mathbf{x}(\boldsymbol{\mu}) ; \boldsymbol{\theta}_{y_I}\right)=\ell(\mathbf{x}, \boldsymbol{\mu}) \hat{y}_{I}\left(\mathbf{x}(\boldsymbol{\mu}) ; \boldsymbol{\theta}_{y_I}\right) + 1$ to naturally satisfy the Dirichlet boundary condition \eqref{eq:test3}, where the length factor function is 
\begin{linenomath*}
\begin{equation*}
    \ell(\mathbf{x}, \boldsymbol{\mu}) = x_0(2 - x_0)x_1(1 - x_1)(\mu_1^2 - (x_0 - 1.5)^2 - (x_1 - 0.5)^2).
\end{equation*}
\end{linenomath*}
The three neural networks $\hat{y}_{I}\left(\mathbf{x}(\boldsymbol{\mu}); \boldsymbol{\theta}_{y_I}\right), \hat{p}_{I}\left(\mathbf{x}(\boldsymbol{\mu}); \boldsymbol{\theta}_{p_I}\right)$ and $\hat{u}\left(\mathbf{x}(\boldsymbol{\mu}); \boldsymbol{\theta}_{u}\right)$ are comprised of three ResNet blocks, each of which contains two fully connected layers with $25$ units and a residual connection. The input dimension of these three neural networks is $4$ and the total number of parameters of these three neural networks is $3 \times 3401 = 10203$. We take $\gamma=0.985$ and the number of epochs for training the state function and the adjoint function increases from $200$ to $700$ during training. The configurations of the neural networks for the PINN and PINN+Projection methods are the same as those of AONN, and the number of training epoch is $50000$. For this test problem, the AONN algorithm converges in $300$ steps. Note that our AONN method can obtain all-at-once solution for any parameter $\boldsymbol{\mu}$. To evaluate the performance of AONN, we employ the classical finite element method to solve the $\mathrm{OCP}(\mb{\mu})$ with a fixed parameter. More specifically, a limited-memory BFGS algorithm implemented with bounded support is adopted in the dolfin-adjoint \cite{mitusch2019dolfin} to solve the corresponding $\mathrm{OCP}$. The solution obtained using the dolfin-adjoint can be regarded as the ground truth. Among the four methods, AONN, PINN and PINN+Projection are able to solve parametric optimal control problems, while the dolfin-adjoint solver can only solve the optimal control problem with a fixed parameter.  

\begin{figure}[htbp]
\centering
\includegraphics[width=14cm]{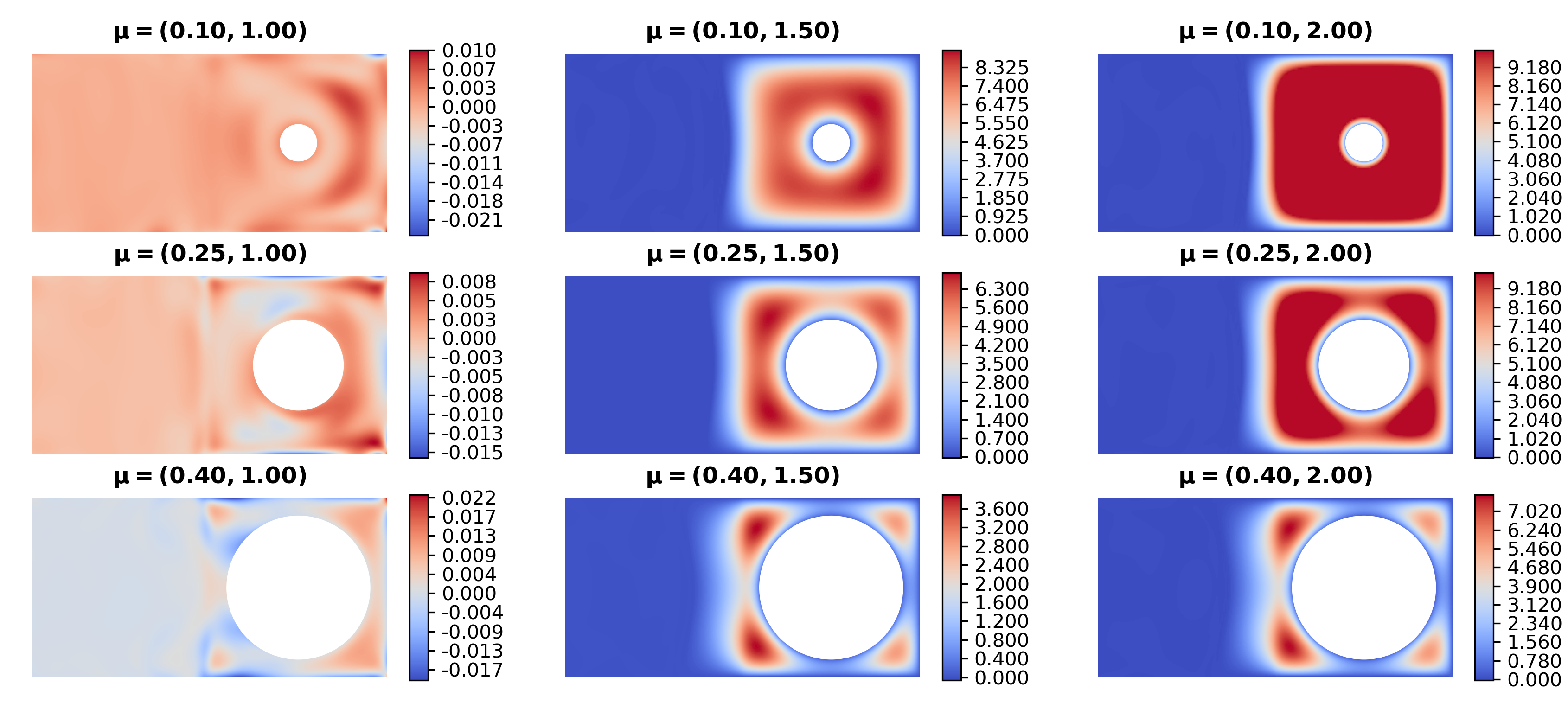}
\caption{Test 4: the AONN solutions $u(\boldsymbol{\mu})$ with several realizations of $\boldsymbol{\mu}=(\mu_1, \mu_2)$.}
\label{fig:test3-2}
\end{figure}

\begin{figure}[htbp]
\centering
\includegraphics[width=14cm]{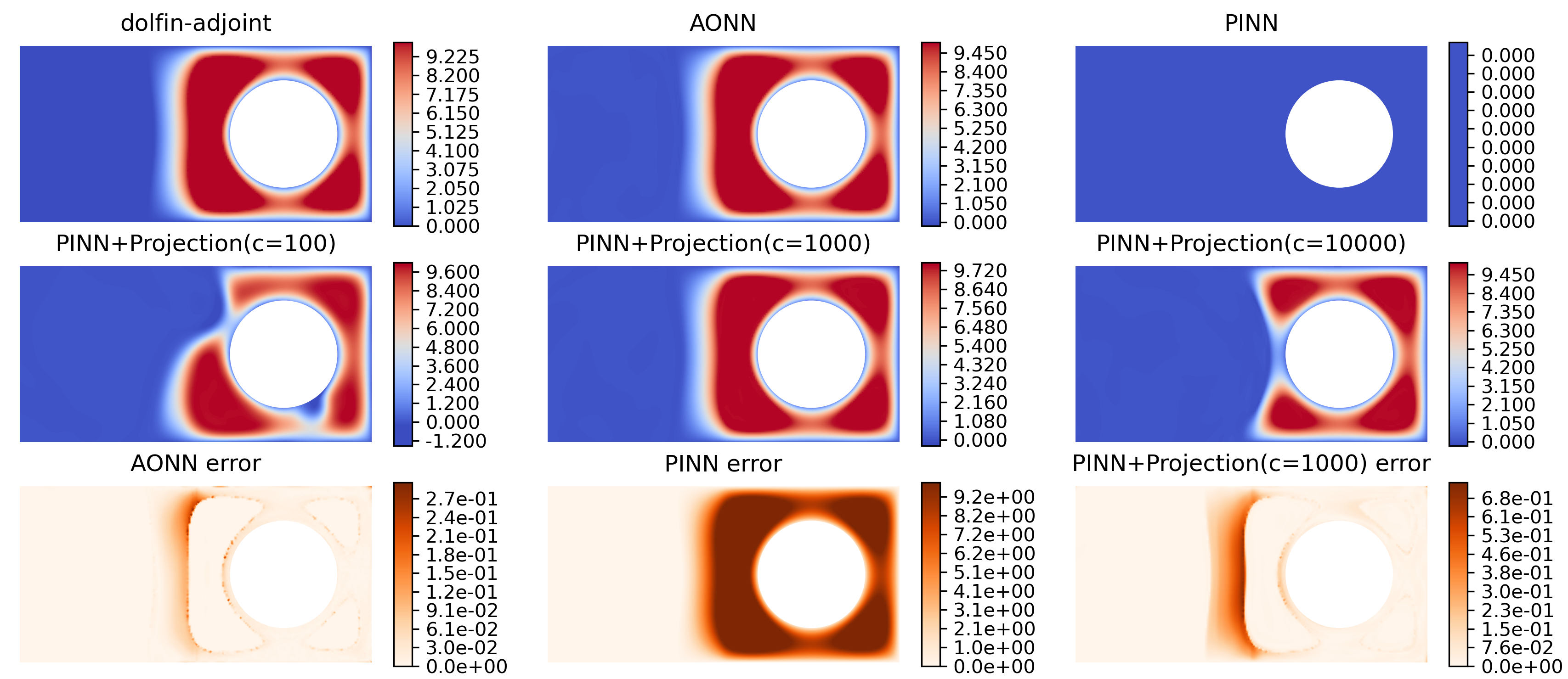}
    \caption{Test 4: the solution obtained by the dolfin-adjoint solver for a fixed parameter $\boldsymbol{\mu}=(0.3, 2.5)$, the approximate solutions of $u$ obtained by AONN, PINN, PINN+Projection (with different $c=100,1000,10000$), and the absolute errors of the AONN solution and the PINN+Projection solution with $c=\frac{1}{\alpha}=1000$.}
\label{fig:test3-3}
\end{figure}

\begin{figure}[!htbp]
  \centering
  \subfigure[]
  {\includegraphics[width=0.31\textwidth]{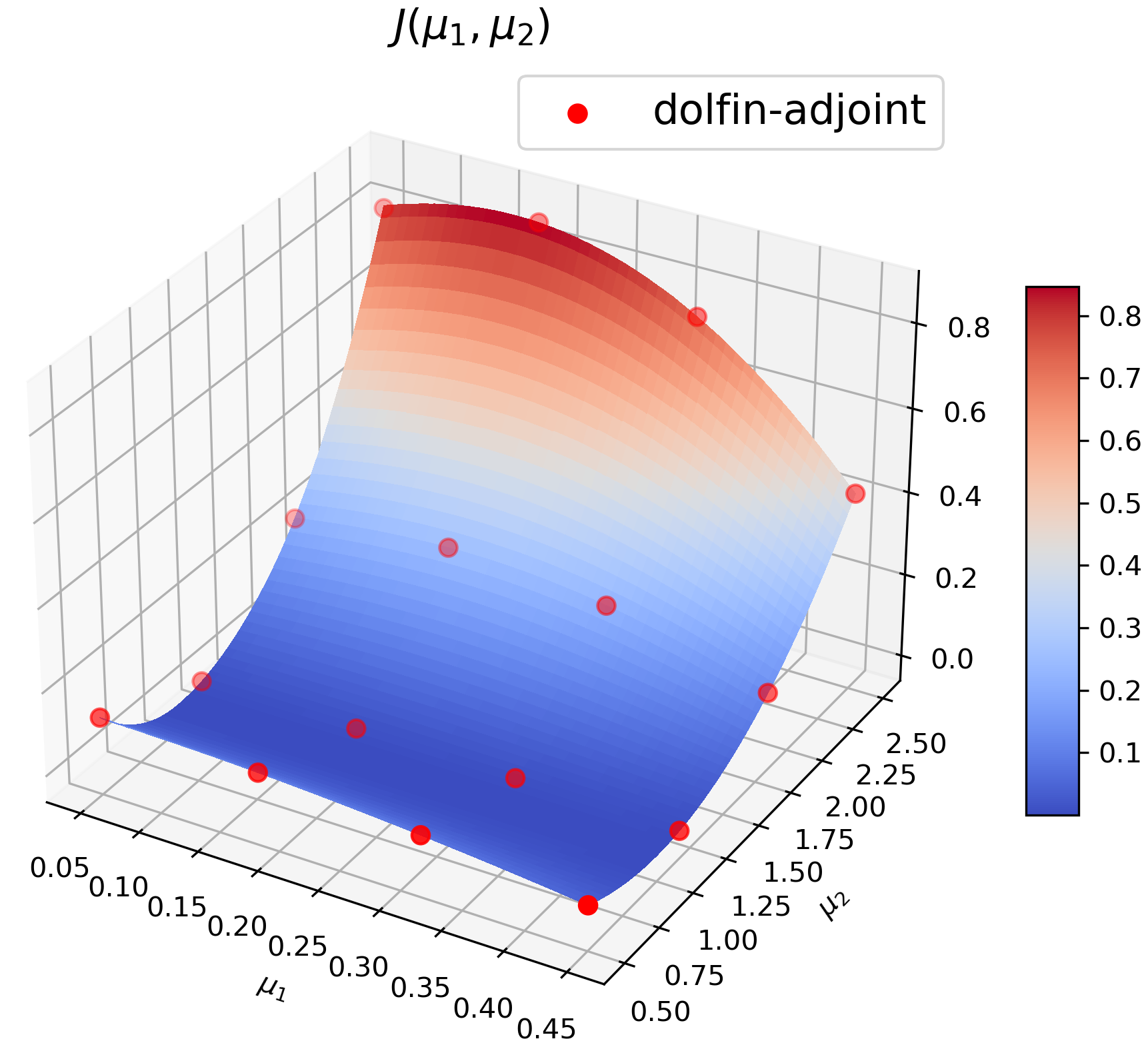}  \label{fig:test3-3d-1}}
  \subfigure[]
  {\includegraphics[width=0.31\textwidth]{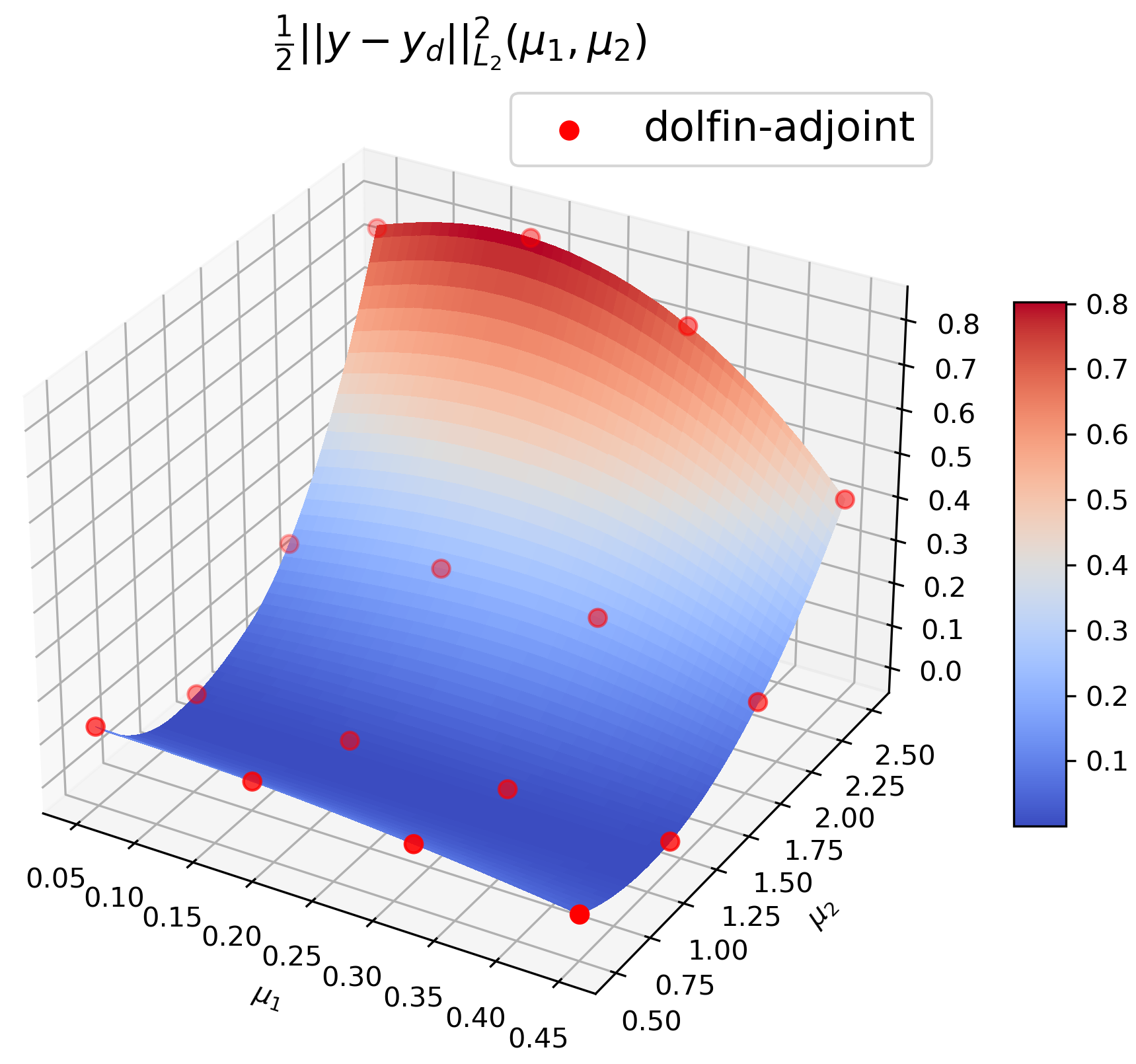}  \label{fig:test3-3d-2}}
  \subfigure[]
  {\includegraphics[width=0.31\textwidth]{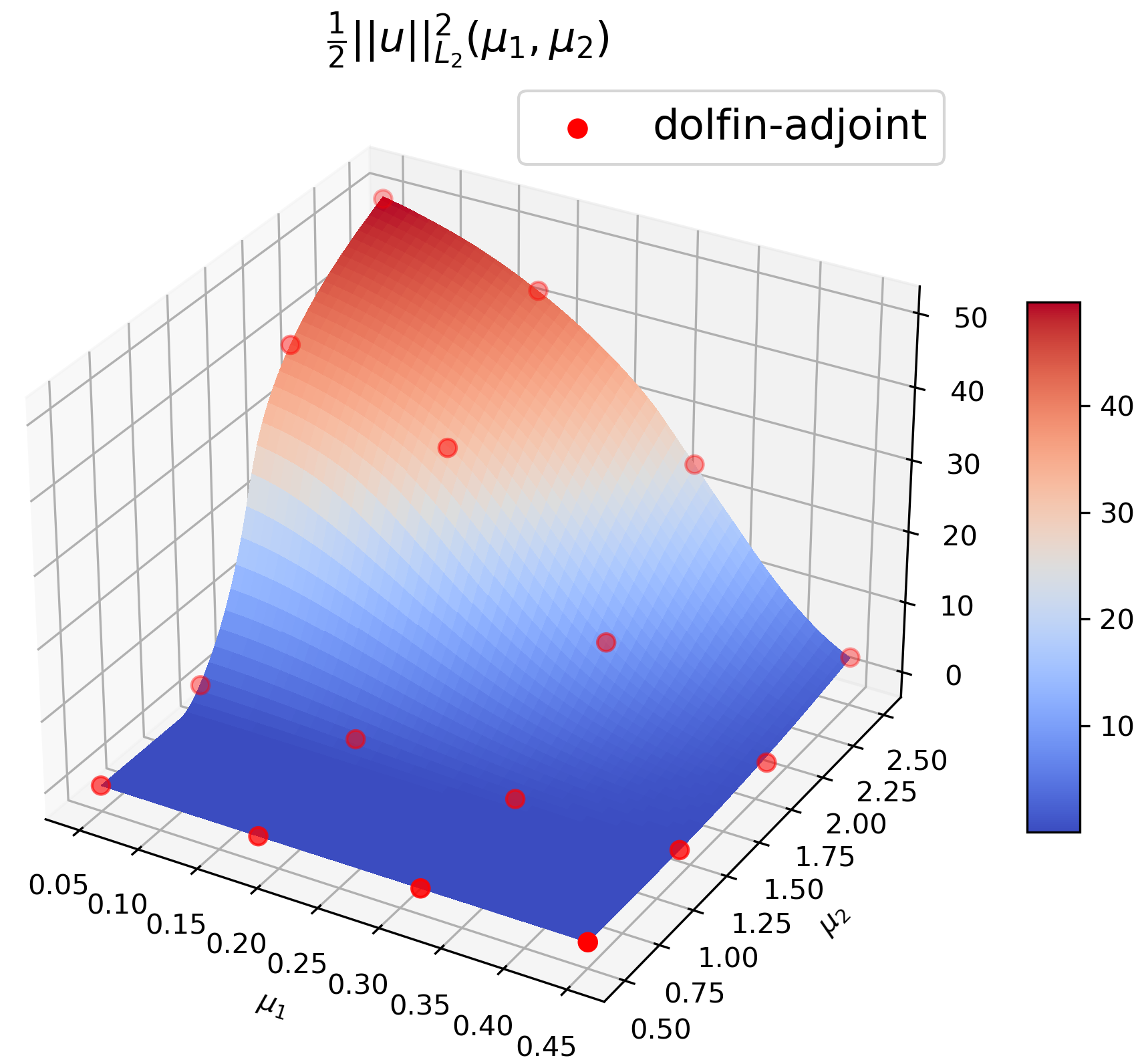}  \label{fig:test3-3d-3}}
  \caption{Test 4: several quantities as functions with respect to parameter $\boldsymbol{\mu}=(\mu_1, \mu_2)$ obtained by AONN. Each red dot denotes the quantity corresponding to a specific $\boldsymbol{\mu}$ computed from the dolfin-adjoint solver. (a) Objective value: $J$ (b) Attainability of the desired state:  $\frac{1}{2}\|y-y_d\|_{L_2}^2$. (c) $L_2$-norm of control function: $\frac{1}{2}\|u\|_{L_2}^2$.}
  \label{fig:test3-3d}
\end{figure}

Figure~\ref{fig:test3-2} shows the optimal control solution obtained using  AONN for the parametric optimal control problem \eqref{eq:test3}. We choose several different parameters $\boldsymbol{\mu}$ for visualization. The left column of Figure~\ref{fig:test3-2} corresponds to $\mu_2=1$, in which case the optimal control is exactly zero because the desired state is achievable for $y=y_d \equiv 1$. The middle and right column of Figure~\ref{fig:test3-2} indicate that the decrease of $\mu_1$ and increase of $\mu_2$ could increase the magnitude of $u$. 
The results obtained by the dolfin-adjoint solver, AONN, PINN and PINN+Projection with different values of $c$ are displayed in Figure~\ref{fig:test3-3}, where the control functions at $\boldsymbol{\mu}=(0.3, 2.5)$ are compared. The mesh with $138604$ triangular elements are used in the dolfin-adjoint solver, and after $16$ steps, the final projected gradient norm achieves $2.379 \times 10^{-10}$. Figure~\ref{fig:test3-3} shows that AONN can converge to the reference solution obtained by the dolfin-adjoint solver but PINN cannot obtain an accurate solution, while the results of PINN+Projection depends heavily on the choice of $c$ in \eqref{projectedKKT}. When $c$ is not equal to $1 / \alpha=1000$, the PINN+Projection method is not guaranteed to converge to the reference solution. This also confirms that the variational loss \eqref{eq_rescontrol} brings great difficulties to neural network training of the KKT system \eqref{projectedKKT}, unless $c=1 / \alpha$, in which case the control function $u$ is canceled out inside the projection operator. 
\begin{linenomath*}
\begin{equation*}
    \mathbf{P}_{U_{ad}(\boldsymbol{\mu})}\left(u(\boldsymbol{\mu})-c \mathrm{d}_{u}J(y(\boldsymbol{\mu}),u(\boldsymbol{\mu});\boldsymbol{\mu})\right) = \mathbf{P}_{U_{ad}(\boldsymbol{\mu})}\left(u(\boldsymbol{\mu})-c(\alpha u(\boldsymbol{\mu})+p(\boldsymbol{\mu}))\right) = \mathbf{P}_{U_{ad}(\boldsymbol{\mu})}\left(-\frac{1}{\alpha}p(\boldsymbol{\mu})\right) .
\end{equation*}
\end{linenomath*}
However, for the next non-smooth test problem, $u$ cannot be separated from the variational loss for any $c$, which results in failure for the PINN+Projection method. 
To demonstrate that AONN can get all-at-once solutions, we first take a $100 \times 100$ grid of $\mathcal{P}$ and choose several different realizations of $\mb{\mu}$ to solve their corresponding OCP using the dolfin-adjoint solver. Then the parameters on the grid together with spatial coordinates are run through the trained neural networks obtained by Algorithm \ref{alg_AONN} to get the optimal solutions of OCP($\mb{\mu}$) all at once. It is worth noting that using the dolfin-adjoint solver to compute the optimal solutions for all parameters on the $100 \times 100$ grid is computationally expensive since $10000$ simulations are required. So we only take $16$ representative points on the grid for the dolfin-adjoint solver (It still takes several hours). Nevertheless, all-at-once solutions can be computed effectively and efficiently through our AONN framework.
Figure~\ref{fig:test3-3d} displays three quantities with respect to $\mu_1, \mu_2$, where Figure~\ref{fig:test3-3d-1} shows the objective functional $J$, Figure~\ref{fig:test3-3d-2} is the accessibility of the desired state and Figure~\ref{fig:test3-3d-3} displays the $L_2$-norm of the optimal control $u$. The red dots in Figure~\ref{fig:test3-3d} show the results obtained by the dolfin-adjoint solver, where $16$ simulations of $\mathrm{OCP}$ with $(\mu_1, \mu_2)\in \{0.05, 0.1833, 0.3167, 0.45\} \times \{0.5, 1.1667, 1.8333, 2.5 \}$ are performed. From Figure~\ref{fig:test3-3d}, it is clear that AONN can obtain accurate solutions.

\subsection{Test 5: Optimal control for the semilinear elliptic equations with sparsity parametrization}\label{sec_test5_sparse}
In this test problem, we again consider a control problem for the semilinear elliptic equations as that in \eqref{test1}. However, this time we consider a sparse optimal control problem with sparsity parametrization. 
The sparse solution in optimal control is often achieved by $L_1$-control cost \cite{casas2011approximation,casas2012optimality,casas2017review} and its application to the controller placement problems is well studied \cite{stadler2009elliptic}. Specifically, we consider the following objective functional with $L_1$-control cost:
\begin{linenomath*}
\begin{equation*}
J(y,u) = \frac{1}{2}\left\|y-y_{d}\right\|_{L_{2}}^{2}+\frac{\alpha}{2}\|u\|_{L_{2}}^{2}+\beta\|u\|_{L_{1}},
\end{equation*}
\end{linenomath*}
where the coefficient $\beta$ of the $L_1$-term controls the sparsity of the control function $u$. With the increase of $\beta$, the optimal control gradually becomes sparse and eventually reaches zero. In order to make continuous observation of this phenomenon, we need to solve the following parametric optimal control problem by setting $\beta$ as a variable parameter $\boldsymbol{\mu}$,
\begin{linenomath*}\begin{equation}
    \left\{\begin{aligned} 
    &\min_{y(\boldsymbol{\mu}),u(\boldsymbol{\mu})}  J(y(\boldsymbol{\mu}), u(\boldsymbol{\mu}); \boldsymbol{\mu}):=\frac{1}{2}\left\|y(\boldsymbol{\mu})-y_{d}\right\|_{L_{2}(\Omega)}^{2}+\frac{\alpha}{2}\|u(\boldsymbol{\mu})\|_{L_{2}(\Omega)}^{2}+\boldsymbol{\mu}\|u(\boldsymbol{\mu})\|_{L_{1}(\Omega)},\\
    &\text { subject to }  
    \left\{\begin{aligned}
    -\Delta y(\boldsymbol{\mu}) + y(\boldsymbol{\mu})^3&=u(\boldsymbol{\mu})&&\text{ in } \Omega, \\ 
    y(\boldsymbol{\mu}) &=0  &&\text { on } \partial \Omega,
    \end{aligned}\right.\\ 
    &\text{and} \quad u_a \leq u(\boldsymbol{\mu}) \leq u_b \quad \text { a.e. in } \Omega.\\
    \end{aligned}\right.
\end{equation}\end{linenomath*}
We fixed the other parameters
\begin{linenomath*}
\begin{equation*}
\begin{aligned}
&\Omega=B(0,1),\\
&\alpha=0.002, u_a=-12, u_b=12,\\
&y_d = 4\sin \left(2\pi x_{1}\right) \sin \left(\pi x_{2}\right)\exp(x_1),\\
\end{aligned}
\end{equation*}
\end{linenomath*}
and the range of parameter is set to $\boldsymbol{\mu} \in [0, \boldsymbol{\mu}_{max}]$. The upper bound $\boldsymbol{\mu}_{max} = 0.128$ ensures that for any $\boldsymbol{\mu} \geq \boldsymbol{\mu}_{max}$ the optimal control $u^*(\boldsymbol{\mu})$ is identically zero. We compute the generalized derivative
 \begin{linenomath*}\begin{equation}\label{subgradient}
    \mathrm{d}_{u}J(y(\boldsymbol{\mu}), u(\boldsymbol{\mu}); \boldsymbol{\mu}) = \alpha u(\boldsymbol{\mu}) + p(\boldsymbol{\mu}) + \boldsymbol{\mu} \ \mathrm{sign}(u(\boldsymbol{\mu})).
\end{equation}\end{linenomath*}
where $p$ is the solution of the adjoint equation as defined in \eqref{test1-adjoint}, and $\mathrm{sign}$ is an element-wise operator that extracts the sign of a function. 

Since the optimal control function varies for different $\boldsymbol{\mu}$, solving a series of sparse optimal control problems is straightforward in general. For example, Eduardo Cases \cite{casas2017review} calculated the optimal solutions for $\boldsymbol{\mu} = 2^i \times 10^{-3}, i=0,1,\ldots,8$. Here, we use AONN to compute all the optimal solutions for any $\boldsymbol{\mu} \in [0,0.128]$ all at once. The length factor function for the Dirichlet boundary condition is chosen as $\ell(\mathbf{x}) = 1 - x_0^2 - x_1^2$. The neural networks 
$\hat{y}_{I}\left(\mathbf{x}(\boldsymbol{\mu}); \boldsymbol{\theta}_{y_I}\right), \hat{p}_{I}\left(\mathbf{x}(\boldsymbol{\mu}); \boldsymbol{\theta}_{p_I}\right)$ and $\hat{u}\left(\mathbf{x}(\boldsymbol{\mu}); \boldsymbol{\theta}_{u}\right)$
are trained by AONN, which have the same configurations to those in the previous test, except for the input dimension being $3$, resulting in $3 \times 3376=10128$ undecided parameters. To this end, we sample $N=20000$ points in the spatio-parametric space $\Omega_{\mathcal{P}} = B(0,1)\times[0,\boldsymbol{\mu}_{max}]$ by a uniform distribution. In order to capture the information at the boundary of $\boldsymbol{\mu}$, $2000$ of these $20000$ points are sampled in $B(0,1)\times\{0\}$ and $B(0,1)\times\{\boldsymbol{\mu}_{max}\}$. We take $500$ iteration steps and gradually increase the training epochs with $n_{\mathrm{aug}}=100$ after every $100$ iteration. As a result, training epochs for the state function and the adjoint function increase from $200$ to $600$ during training. The step size $c^k$ starts with $c^0 = 10$ and decreases by a factor $\gamma = 0.985$ after every iteration.

\begin{figure}[htbp]
\centering
\includegraphics[height=10cm]{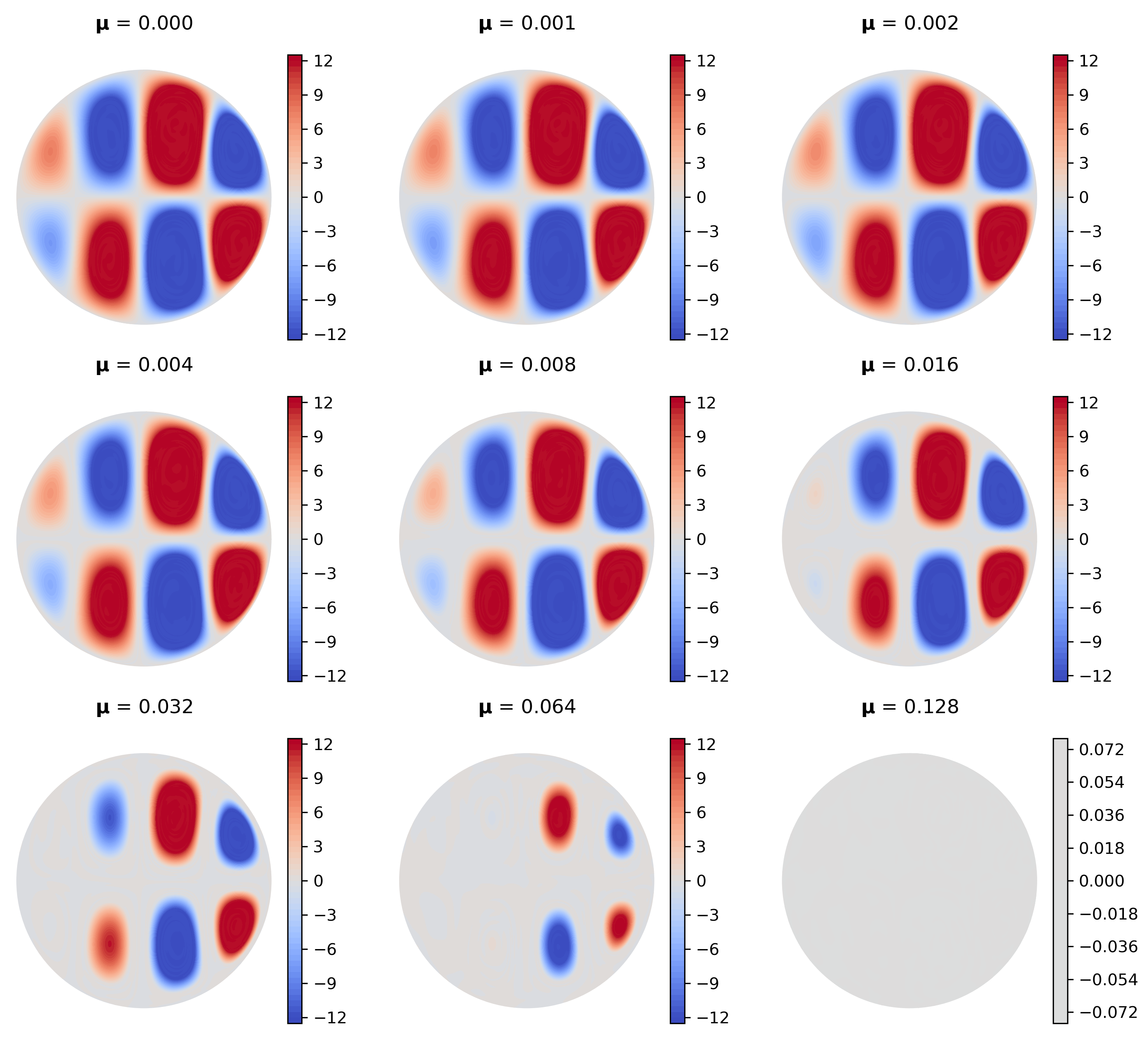}
\caption{Test 5: the AONN solutions $u(\boldsymbol{\mu})$ of representative values for $\boldsymbol{\mu} = 2^i \times 10^{-3}, i=0,1,\ldots,8$.}
\label{fig:test2-1}
\end{figure}

\begin{figure}[htbp]
\centering
\includegraphics[height=6cm]{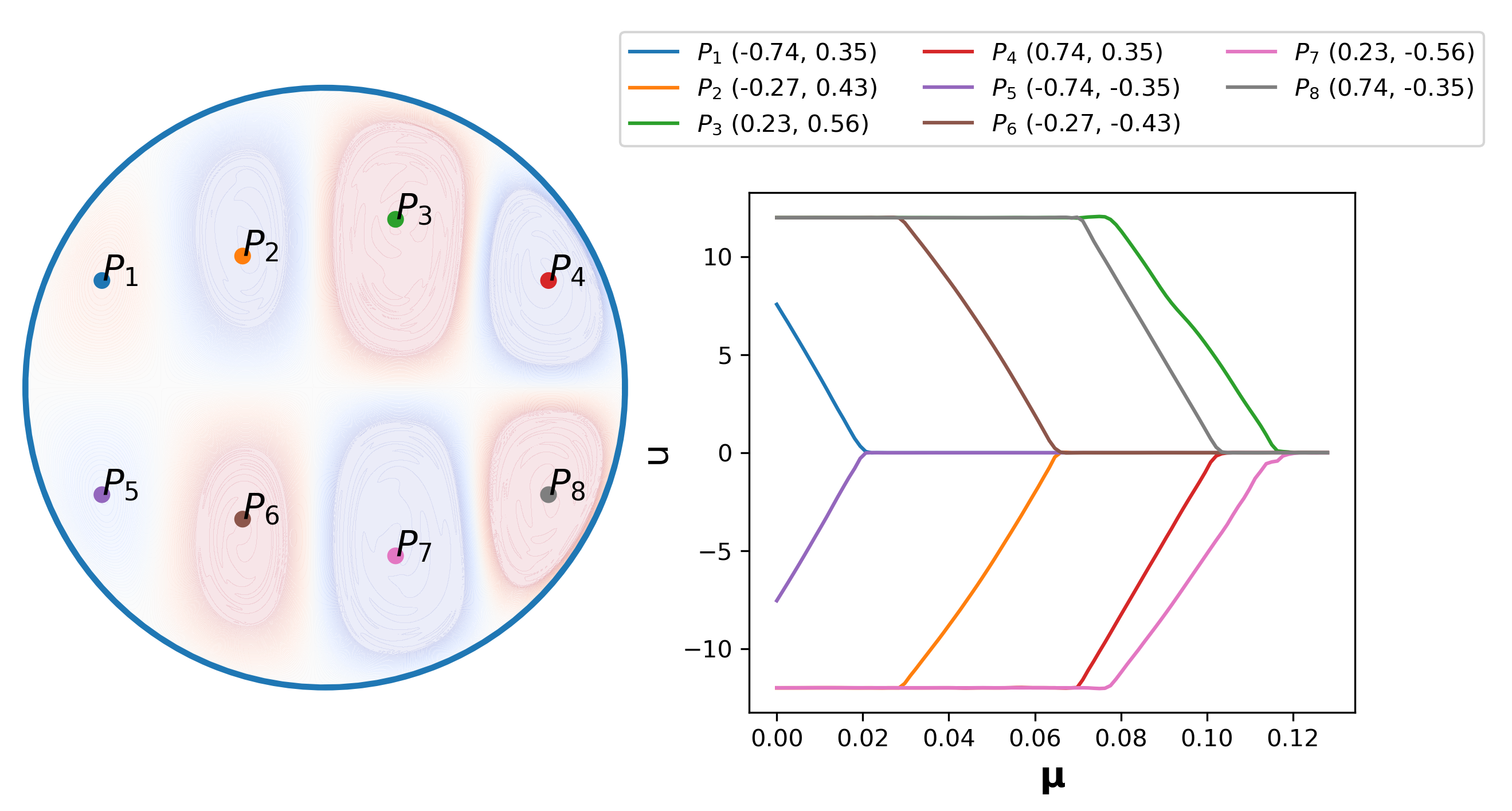}
\caption{Test 5: the AONN solution $u(\boldsymbol{\mu})$ of eight fixed peaks $P_1\sim P_8$ as a function respect to $\boldsymbol{\mu}$. The legend on the right is the coordinates of the eight points.}
\label{fig:test2-2}
\end{figure}

The optimal control for some representative $\boldsymbol{\mu}\in[0,\boldsymbol{\mu}_{max}]$ computed by AONN are displayed in Figure~\ref{fig:test2-1}. The AONN results are consistent with the results presented in ref.\cite{casas2017review}, where the sparsity of optimal control increases as $\boldsymbol{\mu}$ increases. As shown in Figure~\ref{fig:test2-1}, the initial optimal control for $\boldsymbol{\mu}=0$ has eight peaks and each peak disappears as $\boldsymbol{\mu}$ increases. 
To determine where it is most efficient to put the control device, one might require some manual tuning of $\boldsymbol{\mu}$ and thus need to solve $\mathrm{OCP}$ many times for different $\boldsymbol{\mu}$.
Determining these optimal locations is easy if we have obtained the parametric solutions $u^*(x,y,\boldsymbol{\mu})$, which is exactly what AONN does. 
The coordinates of the eight peaks are obtained by evaluating the last vanishing positions of $u^*(x,y,\boldsymbol{\mu})$ as $\boldsymbol{\mu}$ increases at a uniform $100^3$ grid on $[-1,1]\times[-1,1]\times[0,\boldsymbol{\mu}_{max}]$.
Figure~\ref{fig:test2-2} shows the variation of control values at the eight peaks as a function with respect to $\boldsymbol{\mu}$. We observe that the control function values at points $P_1,P_2,P_3$ start with $12$, and begin to decrease after $\boldsymbol{\mu}$ reaches a certain value, finally drop to zero. Value at $P_4$ starts to decrease from a number less than $12$ until it reaches zero. The behavior of points $P_5 \sim P_8$ is completely symmetric.
 
To conclude, with these five numerical tests, we examine the efficiency of AONN and compare its performance with PINN, PINN+Projection, and the traditional solver. The numerical results indicate that the proposed AONN method is more advantageous than the PINN+Projection method and the PINN method in solving parametric optimal control problems. The PINN method cannot obtain accurate solutions for complex constrained problems, and the PINN+projection method improves the accuracy of the PINN method in general but has limitations on nonsmooth problems such as the sparse optimal control problems, while the AONN method is a general framework performing better on different types of parametric optimal control problems.
\section{Conclusions}
\label{sec:conclusions}
We have developed AONN, an adjoint-oriented neural network method, for computing all-at-once solutions to parametric optimal control problems. That is, the optimal control solutions for arbitrary parameters can be obtained by solving only once. 
The key idea of AONN is to employ three neural networks to approximate the control function, the adjoint function, and the state function in the optimality conditions, which allows this method to integrate the idea of the direct-adjoint looping (DAL) approach in neural network approximation. In this way, three parametric surrogate models using neural networks provide all-at-once representations of optimal solutions, which avoids mesh generation for both spatial and parametric spaces and thus can be generalized to high-dimensional problems. With the integration of DAL, AONN also avoids the penalty-based loss function of the complex Karush–Kuhn–Tucker (KKT) system, thereby reducing the training difficulty of neural networks and improving the accuracy of solutions. 
Numerical experiments have shown that AONN can solve parametric optimal control problems all at once with high accuracy in several application scenarios, including control parameters, physical parameters, model parameters, and geometrical parameters. 

Many questions remain open, e.g., choosing the step size and the scaling factor are heuristic, and solving some complex problems requires a higher computational cost.
Future works could include the analysis of the convergence rate to better understand the properties of AONN, the introduction of adaptive sampling strategies to further improve both robustness and effectiveness, and the generalization and application of AONN to more challenging problems such as shape or topology optimizations.

\section*{Acknowledgments}

This study was funded in part by National Natural Science Foundation of China (\#12131002) and China Postdoctoral Science Foundation (2022M711730).

\bibliographystyle{siamplain}
\bibliography{references}

\end{document}